\documentclass{amsart}   
\usepackage{amsmath}
\usepackage[mathscr]{eucal} 
\usepackage{amssymb}
\usepackage{latexsym}
\usepackage{amsthm} 
\theoremstyle{plain}
\newtheorem{theorem}{Theorem}

\newtheorem{proposition}{Proposition}

\newtheorem{lemma}{Lemma} 
\newtheorem{assertion}{Assertion $A(n)$\!\!}
\newcommand{\supp}{\mathop{\mathrm{supp}}\nolimits}

\numberwithin{equation}{section}  
\theoremstyle{definition}

\theoremstyle{remark}
\newtheorem{remark}{Remark}

\title[singular integrals on homogeneous groups]
{Estimates for singular integrals on homogeneous groups}  
 \author{Shuichi Sato} 
\begin{document}
\setcounter{page}{1}
\address{Department of Mathematics, 
Faculty of Education, 
Kanazawa University,    
Kanazawa 920-1192, 
Japan} 
\email{shuichi@kenroku.kanazawa-u.ac.jp} 
\thanks{
\endgraf 
}

\maketitle 

\begin{abstract} 
We consider singular integral operators and maximal singular integral operators  with rough kernels on homogeneous groups.  We prove 
certain  estimates for the operators  that imply $L^p$ boundedness of them 
by an extrapolation argument 
under a sharp condition for the kernels. Also, we prove some 
weighted $L^p$ inequalities for the  operators. 
\end{abstract}   
 
\section{ Introduction}   
Let $\Bbb R^n$, $n\geq 2$, be the $n$ dimensional Euclidean space.  
We also regard $\Bbb R^n$ as a homogeneous group with multiplication given 
by a polynomial mapping.  So, we have a dilation family 
$\{A_t\}_{t>0}$ on $\Bbb R^n$ such that each $A_t$ is an automorphism of 
the group structure, where $A_t$ is of the form  
 $$A_tx=(t^{a_1}x_1,t^{a_2}x_2,\dots, t^{a_n}x_n), \quad 
 x=(x_1,\dots, x_n), $$ 
with some real numbers $a_1, \dots , a_n$ satisfying 
 $0<a_1\leq a_2\leq \dots \leq a_n$ 
 (see \cite{T} and \cite[Section 2 of Chapter 1]{NS}).     
 We also write $\Bbb R^n=\Bbb H$.  
In addition to the Euclidean structure, 
$\Bbb H$ is equipped with 
a homogeneous nilpotent Lie group structure, where  Lebesgue measure is a 
bi-invariant Haar measure, the identity is the origin $0$, $x^{-1}=-x$ and 
multiplication $xy$, $x, y \in \Bbb H$, satisfies 
\begin{enumerate}  
\renewcommand{\labelenumi}{(\arabic{enumi})}  
\item (ux)(vx)= ux+vx, $x\in \Bbb H$, $u, v \in \Bbb R$; 
\item $A_t(xy)=(A_tx)(A_ty)$,  $x, y \in \Bbb H$, $t>0$;  
\item if $z=xy$, then $z_k=P_k(x,y)$, where $P_1(x,y)=x_1+y_1$ and 
$P_k(x,y)=x_k+y_k+R_k(x,y)$ for $k\geq 2$ with  a polynomial   
$R_k(x,y)$ depending only on 
$x_1, \dots ,x_{k-1}, y_1, \dots ,y_{k-1}$.   
 \end{enumerate}  
 We denote by $|x|$ the Euclidean norm  for $x\in \Bbb R^n$. Also,   
  we have a norm function $r(x)$  satisfying $r(A_tx)=tr(x)$ 
for  $t>0$ and $x\in \Bbb R^n$. We assume the following:  
\begin{enumerate}  
\renewcommand{\labelenumi}{(\arabic{enumi})}  
\item[(4)]  the function $r$ 
 is continuous on $\Bbb  R^n$ and smooth in 
 $\Bbb R^n\setminus \{0\}$;  
\item[(5)]  $r(x+y)\leq C_0(r(x)+r(y))$, $r(xy)\leq C_0(r(x)+r(y))$ for some 
 constant $C_0\geq 1$, $r(x^{-1})=r(x)$;  
\item[(6)]  
there are positive constants $c_1, c_2, c_3, c_4, \alpha_1, \alpha_2, 
\beta_1$ and $\beta_2$ such that 
\begin{gather*} 
c_1|x|^{\alpha_1}\leq r(x)\leq c_2|x|^{\alpha_2} \quad \text{if $r(x)\geq 1$}, 
\\ 
c_3|x|^{\beta_1}\leq r(x)\leq c_4|x|^{\beta_2} 
\quad \text{if $r(x)\leq  1$};  
\end{gather*}  
\item[(7)]  if $\Sigma=\{x\in \Bbb R^n: r(x)=1\}$, then 
$\Sigma=S^{n-1}=\{x\in \Bbb R^n: |x|=1\}$.    
 \end{enumerate}    
 Let $\gamma=a_1+\dots +a_n$. 
 Then,  $dx=t^{\gamma-1}\ dS\,dt$,  that is,  
$$\int_{\Bbb R^n}f(x)\,dx=\int_0^\infty\int_\Sigma f(A_t\theta)t^{\gamma-1}\,
dS(\theta)\,dt $$ 
for an appropriate function $f$ with $dS=\omega\, dS_0$, where $\omega$ is a 
strictly positive $C^\infty$ function on $\Sigma$ and $dS_0$ is the Lebesgue 
surface measure on $\Sigma$.    
For appropriate functions $f, g$ on $\Bbb H$, the convolution $f*g$ is defined 
by 
$$f*g(x)=\int_{\Bbb R^n} f(y)g(y^{-1}x)\, dy. $$  
The space $\Bbb H$ with a left invariant quasi-metric $d(x,y)=r(x^{-1}y)$ 
can be regarded as a space of homogeneous type 
(see \cite{CT, CW, FR, NS, R, St, SW, T} for more details). 
 \par   
The Heisenberg group $\Bbb H_1$ is an example of a homogeneous group. 
If we define the multiplication 
$$(x,y,u)(x',y',u')=(x+x',y+y',u+u'+(xy'-yx')/2), $$  
 $(x,y,u), (x',y',u')\in \Bbb R^3$, then $\Bbb R^3$ with this group law 
 is the  Heisenberg group $\Bbb H_1$; a dilation is defined by 
 $A_t(x,y,u)=(tx,ty,t^2u)$. 
 \par 
Let $\Omega$ be locally integrable in $\Bbb R^n\setminus\{0\}$ and 
homogeneous of degree $0$ with respect to the dilation group 
$\{A_t\}$, that is,  $\Omega(A_tx)=\Omega(x)$ for $x\neq 0$, $t>0$. 
We assume  that 
$$ \int_\Sigma \Omega(\theta)\, dS(\theta)=0.$$ 
Let $K(x)=\Omega(x')r(x)^{-\gamma}$,  $x'=A_{r(x)^{-1}}x$ for $x\neq 0$. 
 For $s\geq 1$, let $d_s$ denote the collection of measurable functions 
$h$ on $\Bbb R_+=\{t\in \Bbb R : t>0\}$ satisfying 
$$\|h\|_{d_s}=
\sup_{j\in \Bbb Z}\left(\int_{2^j}^{2^{j+1}}|h(t)|^s\,dt/t\right)^{1/s}
<\infty,$$  
where $\Bbb Z$ denotes the set of integers. 
We define $\|h\|_{d_\infty}=\|h\|_{L^\infty(\Bbb R_+)}$.   
Note that $d_s\subset d_u$ if $s\geq u$.    
Also,  put for  $t\in (0,1]$,    @
$$\omega(h,t)=
 \sup_{|s|<tR/2}\int_R^{2R}|h(r-s)-h(r)|\,dr/r,$$  
where the supremum is taken over all $s$ and $R$ such 
that $|s|<tR/2$ (see \cite{FS, Se}).  
For $\eta>0$, let $\Lambda^\eta$ denote the family of functions  $h$ 
such that 
$$\|h\|_{\Lambda^\eta}=\sup_{t\in(0,1]}t^{-\eta}\omega(h,t)<\infty. $$  
Define a space $\Lambda_s^\eta=d_s\cap \Lambda^\eta$ and set 
$\|h\|_{\Lambda_s^\eta}=\|h\|_{d_s}+\|h\|_{\Lambda^\eta}$ for 
$h\in \Lambda_s^\eta$.  Note that 
$\Lambda_s^{\eta_1}\subset  \Lambda_s^{\eta_2}$ if $\eta_2\leq \eta_1$, and 
$\Lambda_{s_1}^{\eta}\subset  \Lambda_{s_2}^{\eta}$ if $s_2\leq s_1$.  
\par 
Let 
\begin{equation}\label{s1} 
  Tf(x)=\mathrm{p. v.}f*L(x)
  =\mathrm{p. v.} \int_{\Bbb R^n}f(y)L(y^{-1}x)\,dy,    
\end{equation} 
where $L(x)=h(r(x))K(x)$, $h \in d_1$.   
We consider  $L^q(\Sigma)$ spaces and write $\|F\|_q=
  \left(\int_{\Sigma}|F(\theta)|^q\,dS(\theta)\right)^{1/q}$ for 
  $F\in  L^q(\Sigma)$ ($\|F\|_\infty$ is defined as usual).  
Let $s'=s/(s-1)$ denote the conjugate exponent to $s$. 
We shall prove $L^p$ estimates for $Tf$ with $h\in \Lambda_s^{\eta/s'}$ 
and $\Omega\in L^s(\Sigma)$, $s>1$,  as $s$ approaches $1$. 
\begin{theorem}   Let $s>1$.  
Suppose that $\Omega\in L^s(\Sigma)$ and  $h\in \Lambda_s^{\eta/s'}$ 
 for some fixed positive number $\eta$.  
Then, if $1<p<\infty$, 
$$\|Tf\|_p\leq C_ps(s-1)^{-1}\|h\|_{\Lambda_s^{\eta/s'}}\|\Omega\|_s\|f\|_p,  
$$  
where the constant $C_p$ is independent of $s$, $\Omega$ and $h$. 
\end{theorem} 
\par 
 We denote by $L\log L(\Sigma)$  the Zygmund class of all those 
 functions  $F$ on $\Sigma$ 
  which satisfy 
  $$\int_{\Sigma}|F(\theta)|\log(2+|F(\theta)|)\, dS(\theta)
  <\infty.$$  
\par 
Let $\Lambda$ denote the collection of functions $h$ on 
$\Bbb R_+$ such that there exist a sequence $\{h_k\}_{k=1}^\infty$ of functions on $\Bbb R_+$ and a sequence $\{a_k\}_{k=1}^\infty$ of non-negative real 
numbers  satisfying 
$h=\sum_{k=1}^\infty a_kh_k$, $h_k\in \Lambda_{1+1/k}^{1/(k+1)}$, 
$\sup_{k\geq 1}\|h_k\|_{\Lambda_{1+1/k}^{1/(k+1)}}\leq 1$ and 
$\sum_{k=1}^\infty ka_k<\infty$.  
\par  
Theorem 1 implies the following result. 
\begin{theorem} 
Let $Tf$ be as in  \eqref{s1}. Suppose that $h\in \Lambda$ 
and  $\Omega\in L\log L(\Sigma)$. 
Then, 
$T$ is bounded on $L^p(\Bbb R^n)$  for all $p\in (1,\infty)$. 
\end{theorem} 
\par 
When $h=1$ (a constant function), this is due to \cite{T}.  
See \cite{CZ, CDF, NRW, RS, RS2, R} for relevant results and also 
\cite{S3, Se, T} for  weak $(1,1)$ boundedness.   
\par 
 We also consider the maximal singular integral operator   
  \begin{equation}\label{s2}
  T_*f(x) 
= \sup_{N, \epsilon > 0}\left|\int_{\epsilon<r(y)<N}
f(x y^{-1})L(y)\, dy\right|.      
\end{equation}  
We shall prove  analogs of Theorems 1 and 2 for the operator $T_*$. 
\begin{theorem} Let a number $s$ and functions $h$, $\Omega$ be as in 
Theorem $1$. 
Then  we have 
 $$\|T_*f\|_{p}\leq C_p s(s-1)^{-1}
 \|h\|_{\Lambda_s^{\eta/s'}}\|\Omega\|_{s}\|f\|_{p}$$ 
  for all $p\in (1,\infty)$, 
  where $C_p$ is independent of $s$, $h$ and $\Omega$. 
\end{theorem}   
By Theorem 3 we have the following result. 
\begin{theorem} 
Suppose that $\Omega\in L\log L(\Sigma)$ and $h \in \Lambda$.  
Let $T_*f$ be defined as in \eqref{s2}  by using the functions 
 $\Omega$ and $h$.  
   Then,  $T_*$ is bounded on $L^p(\Bbb R^n)$ for $p\in (1,\infty)$, 
\end{theorem} 
This seems to be novel even in the case when $h=1$.  
If $h=1$,  Theorem 2 can be  proved  by 
interpolation  between $L^2$ estimates and weak $(1,1)$ estimates, both 
of which are given in \cite{T}. 
 For $T_*$ with  $\Omega\in L\log L$, weak $(1,1)$ boundedness is yet 
 to be proved  even in the case $h=1$. 
 \par 
In this note we shall show that results of Tao \cite{T} can be used to obtain 
an analog of a theory of Duoandikoetxea and Rubio de 
Francia \cite{DR} for homogeneous groups which can prove Theorems 1 and 3. 
In our situation, Littlewood-Paley theory (see Lemma 6 in Section 4) 
and interpolation arguments are available as in \cite{DR}, although we cannot 
apply Fourier transform estimates as effectively as in \cite{DR}. We shall show that $L^2$ estimates of Lemma 1 in Section 3 can be used as a substitute for  
Fourier transform estimates if we apply  Cotlar's lemma instead of Plancherel's theorem.  Our methods may extend to the study of some other interesting 
operators in harmonic analysis (see \cite{Ch}, \cite{DR}).  
\par 
 Let $\{B_t\}_{t>0}$, $B_t=t^P=\exp((\log t) P)$, be a dilation group on 
 $\Bbb R^n$, where  
$P$ is an $n\times n$ real matrix whose eigenvalues have positive real 
parts.   Let $N$ be a locally integrable function on $\Bbb R^n\setminus\{0\}$ 
such that $N(B_tx)=t^{-\gamma}N(x)$, $\gamma =\text{{\rm trace} $P$}$,  for
 $t>0$ and $x\in \Bbb R^n\setminus \{0\}$.  
 Let $J(x)=h(r(x))N(x)$ with an appropriate norm function $r(x)$ for 
$\{B_t\}_{t>0}$.  If we define 
$$Sf(x)=\mathrm{p. v.} \int_{\Bbb R^n}f(y)J(x-y)\,dy, $$ 
using Euclidean convolution,  assuming an appropriate cancellation condition 
for $J$, then we can apply  methods of Duoandikoetxea and Rubio de Francia  
\cite{DR} via  Fourier transform estimates to prove $L^p$ boundedness,  
$p\in (1,\infty)$, of $S$ under an  $L\log L$ condition on $\{r(x)=1\}$ 
for $N$ and the condition 
$$ \sup_{j \in \Bbb Z} 
\int_{2^j}^{2^{j+1}}|h(r)|\left(\log(2+|h(r)|)\right)^a\, dr/r <\infty   
$$  
for $h$ with some $a>2$.  Also, a similar result for  maximal singular 
integrals  holds (see \cite{S, S2}).  
\par 
We can also prove some weighted norm estimates for $T$ and $T_*$.  
Let $B$ be a subset of $\Bbb H$ such that 
$$B=\{x\in \Bbb H: r(a^{-1}x)<s\}$$ 
for some  $a\in \Bbb H$ and $s>0$.  Then we call $B$ a ball in $\Bbb H$
 with center $a$ and radius $s$ and  write $B=B(a,s)$.  Note that 
 $|B(a,s)|=cs^\gamma$ with 
$c=|B(0,1)|$, where $|S|$ denotes the Lebesgue measure of a set $S$. 
 Let  $\mathscr{A}_p$, $1<p< \infty$,  be the weight 
 class of Muckenhoupt on $\Bbb H$  defined to be the collection of 
 all  weight functions $w$ on  $\Bbb H$  satisfying  
 $$\sup_B \left(|B|^{-1} \int_B w(x)\,dx\right)\left(|B|^{-1} \int_B
w(x)^{-1/(p-1)}dx\right)^{p-1} < \infty, $$
where the supremum is taken over all balls $B$ in $\Bbb H$ 
(see \cite{C, GGKK}).  Also, the class $\mathscr{A}_1$ is defined 
to be the family of  all  weight functions $w$ on  $\Bbb H$  satisfying 
the pointwise inequality $Mw\leq Cw$ almost everywhere,   
where $M$ denotes the Hardy-Littlewood maximal operator 
$$Mf(x)=\sup_{x\in B}|B|^{-1}\int_B|f(y)|\,dy;  $$ 
the supremum is taken over all balls $B$ in $\Bbb H$ containing  $x$ 
(see \cite{C, CW, GGKK}). 
We can prove the following weighted estimates.   
\begin{theorem}  Let $q>1$.  
Suppose that $\Omega\in L^q(\Sigma) $ and $h\in \Lambda_q^\eta$ for some 
 $\eta>0$.  Let $1<p<\infty$. 
Then, 
\begin {enumerate}
\renewcommand{\labelenumi}{(\arabic{enumi})}   
\item $T$ and $T_*$ are  bounded on $L^p(w)$ if $q' \leq p < \infty$ and 
$w \in \mathscr{A}_{p/q'};$   
\item if $1 < p \leq q$ and $w \in \mathscr{A}_{p'/q'}$, $T$ and  $T_*$ are 
bounded on $L^p(w^{1-p})$.  
\end{enumerate} 
\end{theorem}
See \cite{D, W} for the case of rough singular integrals 
defined by  Euclidean convolution.  
\par 
In Section 2, we shall give some preliminary results from \cite{T} for 
calculation on homogeneous groups. 
A basic $L^2$ estimates (Lemma 1) will be proved in Section 3 by applying 
 methods of  \cite{T}. 
Using the $L^2$ estimate, we shall prove Theorem 1 in Section 4  
by means of a process of \cite{DR, S, S2}.  
In Section 5, we shall prove 
Theorem 3 by adapting arguments of \cite{DR} for the present situation. 
Theorem 5 will be proved in Section 6 by applying arguments of \cite{D} and 
using results of Sections 3--5. 
Finally, we shall prove Theorem 2 from Theorem 1 in Section 7 by 
an extrapolation argument.   
Theorem 4  can be proved in the same way from Theorem 3. 
In what follows, even when we consider functions that may assume general 
complex values,  we deal with real valued functions only to simplify our 
arguments.  
 The letters $C, c$ will be used to denote positive constants which 
 may be different in different occurrences.  

\section{Preliminary results}  
In this section we recall several results from \cite{T}.  Let $f:\Bbb R\to 
\Bbb H$ be smooth. Then the Euclidean derivative 
$\partial_tf(t)$ is defined by 
\begin{equation*}
 f(t+\epsilon)=f(t)+\epsilon\partial_tf(t)+\epsilon^2O(1) \quad \text{for 
 $\epsilon \in (0,1]$.}  
\end{equation*} 
We define the left invariant derivative 
$\partial^L_tf(t)$ by 
\begin{equation*}
 f(t+\epsilon)=f(t)(\epsilon\partial^L_tf(t))+\epsilon^2O(1) \quad \text{for 
 $\epsilon\in (0,1]$.}  
\end{equation*} 
Fix $x\in \Bbb H$ and consider $G_x:\Bbb R^n \to \Bbb R^n$ defined by 
$G_x(y)=xy$. 
Let $JG_x(y)$ be the Jacobian matrix of $G_x$ at $y$. Then  $JG_x(y)$ is 
a lower triangular matrix. The components of $JG_x(y)$ are polynomials in 
$x, y$ and each diagonal component is equal to $1$ (see (3) in Section 1).  
We can see that $\partial^L_tf(t)=
JG_{f(t)^{-1}}(f(t))\partial_tf(t)$. 
\par 
We have the product rule 
\begin{equation}\label{1.1}
\partial^L_t(f(t)g(t))=\partial^L_tg(t)+C[g(t)]\partial^L_tf(t),   
\end{equation}  
where $C[x]: \Bbb R^n\to \Bbb R^n$ is a linear mapping defined by 
\begin{equation*}
x^{-1}(\epsilon v)x=\epsilon C[x]v +\epsilon^2O(1)\quad 
\text{for  $\epsilon\in (0,1]$.} 
\end{equation*} 
We note that  
\begin{equation}\label{1.2} 
\begin{split} 
C[x^{-1}]=&C[x]^{-1}, \quad C[A_tx](A_tv)=A_t(C[x]v), 
\\ 
&|C[x]v|\sim |v| \quad \text{if $|x|\leq 1$}.  
\end{split}  
\end{equation} 
\par 
Define a polynomial mapping $X:\Bbb R^n\to \Bbb R^n$ by 
\begin{equation*}
A_{1+\epsilon}x=x(\epsilon X(x))+\epsilon^2O(1) \quad \text{for 
 $\epsilon\in (0,1]$.}  
\end{equation*}  
Then 
\begin{equation}\label{1.3}
X(A_tx)=A_tX(x), \quad r(X(x))\sim r(x)
\end{equation}  
and 
\begin{equation}\label{1.4}
\partial^L_t(A_{s(t)}f(t))=A_{s(t)}\partial^L_tf(t)+ s'(t)s(t)^{-1}\left(
A_{s(t)}X(f(t))\right),  
\end{equation}  
where $s(t)$ is a strictly positive, smooth function on $\Bbb R_+$.  
Also, $X$ is  a diffeomorphism with Jacobian comparable to $1$.  

\section{ $L^2$ estimates}  
Let $\phi$ be a $C^\infty$ function with compact support in $B(0,1)\setminus 
B(0,1/2)$ satisfying   
$\int \phi =1$,  $\phi(x)=\tilde{\phi}(x)$, $\phi(x)\geq 0$ for all $x\in \Bbb R^n$, where $\tilde{\phi}(x)=\phi(x^{-1})$.  
Define 
$$\Delta_k=\delta_{\rho^{k-1}}\phi-\delta_{\rho^{k}}\phi, \quad k\in\Bbb Z, $$ 
where $\delta_t\phi(x)=t^{-\gamma}\phi(A_t^{-1}x)$ and  $\rho\geq 2$.  
Note that $\supp(\Delta_k)\subset B(0,\rho^k)\setminus B(0,\rho^{k-1}/2)$, 
$ \Delta_k=\tilde{\Delta}_k$ and  $\sum_k\Delta_k=\delta$, where $\delta$ is the delta function.  
Let $\psi_j\in C_0^\infty(\Bbb R)$, $j\in \Bbb Z$,  be such that 
\begin{gather*}
\supp(\psi_j)
\subset \{t\in \Bbb R: \rho^{j}\leq t \leq \rho^{j+2}\},  \quad \psi_j\geq 0, 
\\ 
\sum_{j\in\Bbb Z}\psi_j(t)=1 \quad \text{for $t\neq 0$}, 
\\ 
\left|(d/dt)^m\psi_j(t)\right|\leq c_m|t|^{-m} \quad 
\text{for $m=0,1,2,\dots$}, 
\end{gather*} 
where $c_m$ is a constant independent of $\rho$ (this is possible since 
$\rho\geq 2$).  
\par 
Let 
$$S_jL(x)=(\log 2)^{-1}h(r(x))\int_0^\infty \psi_j(t)\delta_tK_0(x)\,dt/t, $$ 
where 
$$ 
K_0(x)=K(x)\chi_{D_0}(x), \quad 
D_0=\{x\in \Bbb R^n: 1\leq r(x)\leq 2\}.  
$$  
Here $\chi_E$ denotes the characteristic function of a set $E$.  
Then $\sum_{j\in \Bbb Z}S_jL=L$. 
Furthermore, let 
\begin{equation}\label{3.a2}
S_j(F,\ell)(x)=(\log 2)^{-1}\ell(r(x))\int_0^\infty\psi_j(t)\delta_tF(x)\, 
dt/t,  
\end{equation} 
where $F\in L^1(\Bbb R^n)$, $\supp(F)\subset D_0$ and $\ell\in d_1$.  
Let $\Phi$ be a non-negative smooth function  such that 
$\int \Phi(x)\,dx=1$, $\Phi(x^{-1})=\Phi(x)$, $\supp(\Phi)\subset B(0,1)$.  
Define  
\begin{equation}\label{epsilon} 
U_\sigma f=U_\sigma(F,\ell)(f)=\sum_j\sigma_j f*\nu_j,  
\end{equation}
where 
\begin{gather*}
\nu_j(x)=\nu_j(F,\ell)(x)=S_j(F,\ell)(x)-\Phi_j(x), 
\\ 
\Phi_j(x)=\Phi_j(F,\ell)(x)=(\int S_j(F,\ell)\, dx)\delta_{\rho^j}\Phi(x), 
\end{gather*} 
and $\sigma=\{\sigma_j\}$ is an arbitrary sequence such that 
 $\sigma_j=1$ or $-1$.  We note that $\int \nu_j(x)\,dx=0$, $S_jL=
 \nu_j(K_0,h)=S_j(K_0,h)$ and $U_\sigma(K_0,h)(f)=Tf$ if $\sigma_j=1$ for 
 all $j$. 
We prove the following $L^2$ estimates.    
\begin{lemma}  
Suppose that $s>1$,  $F\in L^s(D_0)$ and $\ell\in \Lambda_s^{\eta/s'}$ 
for some fixed positive number $\eta$, where we write $F\in L^s(D_0)$ if 
$F\in L^s(\Bbb R^n)$ and $\supp(F)\subset D_0$. 
 Let $\nu_j=\nu_j(F,\ell)$. 
Then,  for $j, k\in \Bbb Z$ we have 
\begin{equation}\label{2.1} 
\|f*\nu_j*\Delta_{k}\|_2\leq 
C(\log\rho)\min(1, \rho^{-\epsilon(|j-k|-c)/s'})
\|\ell\|_{\Lambda_s^{\eta/s'}}\|F\|_s\|f\|_2 
\end{equation} 
for some positive constants $C, \epsilon$ and $c$ independent of $\rho$, 
$s$,  $\ell$ and $F$.  
\end{lemma}    
\begin{proof} It suffices to prove Lemma $1$ with $\nu$ in place of $\nu_j$, 
assuming $j=0$ on the right hand side of \eqref{2.1}, where   
$$\nu(x)=(\log 2)^{-1}\ell(\rho^jr(x))
\int_0^\infty \psi_j(\rho^jt)\delta_tF(x)\,dt/t
-(\int S_j(F,\ell)\, dx)\Phi(x).$$  
This can be seen from change of variables and 
the formulas: $\delta_t(f*g)=(\delta_tf)*(\delta_tg)$, 
$\delta_{\rho^{-j}}\nu_j=\nu$, $\delta_{\rho^{-j}}\Delta_k=\Delta_{k-j}$.  
\par 
If $k\geq 0$, then from the cancellation condition for $\nu$ and the 
smoothness of $\Delta_k$ we have 
\begin{equation}\label{2.2} 
\|\nu*\Delta_k\|_1\leq C(\log\rho)
\min\left(1,\rho^{-\epsilon k+\tau}\right)\|\ell\|_{d_1}\|F\|_1  
\end{equation} 
for some $\epsilon, \tau>0$, 
which implies the conclusion by Young's inequality, if the constant $c$ is 
large enough. 
\par 
The following result is useful. 
\begin{lemma} 
Suppose that $\ell \in d_q$, $F\in L^q(D_0)$ for some $q\geq 1$. 
Put $S=\delta_{\rho^{-j}}S_j(F,\ell)$. 
Then  
$$\|S\|_q\leq C(\log\rho)\|\ell\|_{d_q}\|F\|_q,  $$ 
where the constant $C$ is independent of $\rho$ and $q$. 
\end{lemma} 
\begin{proof} 
Suppose that $q<\infty$. 
Since $\int_0^\infty\psi_j(\rho^jt)\, dt/t\leq 2\log\rho$, 
H\"{o}lder's inequality implies 
\begin{align*}
&\|S\|_q^q\leq (\log 2)^{-q}(2\log\rho)^{q/q'}
\iint_0^\infty |\ell(\rho^jr(x))|^q \psi_j(\rho^jt)|\delta_tF(x)|^q\, dt/t\, 
dx 
\\ 
&=(\log 2)^{-q}(2\log\rho)^{q/q'}
\iint_0^\infty |\ell(\rho^jtr(x))|^q \psi_j(\rho^jt)|F(x)|^qt^{\gamma(1-q)}
\, dt/t\, dx.    
\end{align*} 
Let $N$ be a positive integer such that $\rho^2\leq 2^{N+1}<2\rho^2$. Then 
\begin{align*}
&\int_0^\infty |\ell(\rho^jtr(x))|^q \psi_j(\rho^jt)t^{\gamma(1-q)}\, dt/t
\\ 
&\leq \sum_{m=0}^N \int_{2^m\rho^jr(x)}^{2^{m+1}\rho^jr(x)} 
|\ell(t)|^q (\rho^{-j}r(x)^{-1}t)^{\gamma(1-q)}\, dt/t.  
\\ 
&\leq \sum_{m=0}^N 2^{m\gamma(1-q)}
\int_{2^m\rho^jr(x)}^{2^{m+1}\rho^jr(x)} |\ell(t)|^q \, dt/t
\\ 
&\leq C(\log\rho)\|\ell\|_{d_q}^q,       
\end{align*} 
where $C\leq 12$.  Collecting results, we get the conclusion for $q<\infty$.  
Also, we easily see that $\|S\|_\infty\leq C\|\ell\|_{d_\infty}\|F\|_\infty$, 
which implies the conclusion for $q=\infty$. 
\end{proof} 
\par 
The estimate \eqref{2.2} can be shown as follows. 
First, by Lemma 2 with $q=1$   
\begin{equation}\label{2.3} 
\|\nu*\Delta_k\|_1\leq C(\log\rho)\|\ell\|_{d_1}\|F\|_1. 
\end{equation} 
Suppose that $k\geq 1$. 
Let $t=\rho^{k-1}$. Then $\Delta_k=\delta_t\Delta_1$.  
Since $\int \nu =0$, 
\begin{align}\label{2.4} 
\nu*\Delta_k(x)&=\int t^{-\gamma}\left(\Delta_1(A_t^{-1}(y^{-1}x))-
\Delta_1(A_t^{-1}x))\right)\nu(y)\,dy 
\\  
&=\int t^{-\gamma}\left(\int_0^1 W'(u) \,du\right) \nu(y)\,dy,  \notag 
\end{align}  
where 
$$W(u)=\Delta_1((uY)^{-1}X)=\Delta_1(P_1(-uY,X),\dots, P_n(-uY,X)),$$ 
with $Y=A_{t^{-1}}y, X=A_{t^{-1}}x$ (see (3) in Section 1). 
Note that 
\begin{multline*}
W'(u)=\langle (\nabla\Delta_1)(P_1(-uY,X),\dots, P_n(-uY,X)), 
\\ 
(\partial uP_1(-uY,X),\dots, \partial uP_n(-uY,X))\rangle, 
\end{multline*} 
where 
$\nabla\Delta_1=(\partial x_1\Delta_1, \partial x_2\Delta_1, \dots, 
\partial x_n\Delta_1)$ and $\langle\cdot,\cdot\rangle$ 
denotes the Euclidean inner product in $\Bbb R^n$.    
Also, note that 
$$\partial uP_i(-uY,X)=\langle -Y, \nabla_xP_i(-uY,X)\rangle, $$ 
where 
$$\nabla_xP_i(x,y)=(\partial x_1P_i(x,y), \dots, \partial x_nP_i(x,y)). $$ 
We may assume that $r(Y)\leq C\rho^2$, 
$r(X)\leq C\rho^2$ in \eqref{2.4} by checking the support condition. Therefore 
$$\sup_{u\in [0,1]}|\partial uP_i(-uY,X)|\leq C|Y|\rho^M$$ 
for some $M>0$ and hence 
$$|W'(u)|\leq C|Y|\rho^M|(\nabla\Delta_1)((uY)^{-1}X)|. $$ 
Note that $\|\nabla \Delta_1\|_1\leq C$, $\|A_t^{-1}\|$ 
(a norm for $A_t^{-1}$ as a linear transformation on $\Bbb R^n$) is less than 
$Ct^{-\beta}$ and $|y|\leq C\rho^{2\alpha}$ on the support of $\nu$,  
with $\beta=1/\beta_1, \alpha=1/\alpha_1$ (see (6) of Section 1). Therefore,  
\begin{align}\label{2.5} 
\|\nu*\Delta_k\|_1&\leq C\rho^M\|\nabla\Delta_1\|_1\int |A_t^{-1}(y)||\nu(y)|
\,dy 
\\  
&\leq C\rho^Mt^{-\beta}\int |y||\nu(y)|\,dy    \notag 
\\  
&\leq C\rho^Mt^{-\beta}\rho^{2\alpha}\|\nu\|_1\leq 
C(\log\rho)\rho^{-k\beta}\rho^{2\alpha+\beta+M}\|\ell\|_{d_1}\|F\|_1.    
\notag 
\end{align}  
By \eqref{2.3} and \eqref{2.5}, we have \eqref{2.2} for $k\geq 0$. 
\par 
We next assume that $k\leq -1$.  
Since $\int \Delta_0(x)\,dx=0$, as in the proof of 
\eqref{2.2} we have 
$$\|\Phi*\Delta_k\|_1=\|\delta_{\rho^{-k}}\Phi*\Delta_0\|_1 
\leq C\rho^{-\epsilon|k|}. $$  
Therefore, separately estimating $\|f*S*\Delta_k\|_2$ and 
$\|f*(\int S_j(F,\ell))\Phi*\Delta_k\|_2$, 
it suffices to prove 
\begin{equation}\label{3.a1}
\|f*S*\Delta_{k}\|_2\leq 
C(\log\rho)\min(1, \rho^{-\epsilon(|k|-c)/s'})
\|\ell\|_{\Lambda_s^{\eta/s'}}\|F\|_s\|f\|_2,   
\end{equation} 
where $S=\delta_{\rho^{-j}}S_j(F,\ell)$ as above. 
\par 
By the estimate 
$$\|S*\Delta_k\|_1\leq C(\log\rho)\|\ell\|_{d_1}\|F\|_1  
 $$
and the $T^*T$ method, to prove \eqref{3.a1} it suffices to show that  
\begin{equation}\label{2.6}
\left\|f*\Delta_k*\tilde{S}*S*\Delta_k\right\|_2
\leq C(\log\rho)^2\rho^{\epsilon(k+c)/s'}\|\ell\|_{\Lambda_s^{\eta/s'}}^2
\|F\|_s^2\|f\|_2.      
\end{equation} 
Since $\|T^*T\|=\|(T^*T)^n\|^{1/n}$, \eqref{2.6} follows from 
\begin{equation*} 
\left\|f*\left(\Delta_k*\tilde{S}*S*\Delta_k\right)_*^n\right\|_2 
\leq C(\log\rho)^{2n}\rho^{\epsilon(k+c)/s'}\|\ell\|_{\Lambda_s^{\eta/s'}}^{2n}
\|F\|_s^{2n}\|f\|_2
\end{equation*} 
for some $\epsilon, c >0$, where 
$\left(\Delta_k*\tilde{S}*S*\Delta_k\right)_*^n$ denotes 
the convolution product of $n$ factors of $\Delta_k*\tilde{S}*S*\Delta_k$. 
By Young's inequality, this follows from the $L^1$ estimate 
\begin{equation}\label{2.7} 
\left\|\left(\Delta_k*\tilde{S}*S*\Delta_k\right)_*^n\right\|_1 
\leq C(\log\rho)^{2n}\rho^{\epsilon(k+c)/s'}
\|\ell\|_{\Lambda_s^{\eta/s'}}^{2n}\|F\|_s^{2n}.   
\end{equation} 
Note that $\left(\Delta_k*\tilde{S}*S*\Delta_k\right)_*^n = 
\Delta_k*\tilde{S}*S*\left(\Delta_k*\Delta_k*\tilde{S}*S\right)_*^{n-1}
*\Delta_k$.   
Since $\|\Delta_k'*\tilde{S}\|_1\leq C(\log\rho)\|\ell\|_{d_1}\|F\|_1$ and 
$\Delta_k'*\tilde{S}(x)=\int\delta_y(x)\Delta_k'*\tilde{S}(y)\,dy$, 
where $\delta_y(x)$ is the delta function concentrated at $y$ and $\Delta_k'$ 
is either $\Delta_k$ or $\Delta_k*\Delta_k$,  
\eqref{2.7} follows from 
\begin{equation*} 
\left\|\delta_{w_1}*S*\dots *\delta_{w_n}*S*\Delta_k\right\|_1 \leq 
C(\log\rho)^{n}\rho^{\epsilon(k+c)/s'}\|\ell\|_{\Lambda_s^{\eta/s'}}^n
\|F\|_s^{n}     
\end{equation*}  
uniformly for $w_1, \dots, w_n\in B(0, C\rho^2)$.  
To get this, it suffices to prove 
\begin{equation} \label{2.8}
\left|\left\langle \delta_{w_1}*S*\dots *\delta_{w_n}*S*\Delta_k, g
\right\rangle\right| 
\leq C(\log\rho)^{n}\rho^{\epsilon(k+c)/s'}\|\ell\|_{\Lambda_s^{\eta/s'}}^n
\|F\|_s^{n} 
\end{equation} 
uniformly in $w_1, \dots, w_n\in B(0, C\rho^2)$, 
for all smooth $g$ with compact support satisfying 
$\|g\|_\infty\leq 1$.   
\par 
Fix $g$. Then, the inner product on the left hand side of \eqref{2.8} is equal 
to 
\begin{equation} \label{2.9}
\iiint \Delta_k(x)g(H(y,t)x)\prod_{i=1}^n\left(\ell(t_i,y_i)F(y_i)
\psi(t_i)\right)\, dy\,\bar{d}t\,dx, 
\end{equation}
where $\ell(t_i,y_i)=\ell(t_i\rho^jr(y_i))$, 
$\psi(t_i)=(\log 2)^{-1}\psi_j(\rho^jt_i)$, 
$t=(t_1,\dots, t_n)$, $y=(y_1,\dots,y_n)\in D_0^n$, 
$\bar{d}t=(dt_1/t_1)\dots(dt_n/t_n)$, $dy=dy_1\dots dy_n$ 
and 
$$H(y,t)=\prod_{i=1}^n w_iA_{t_i}y_i=w_1A_{t_1}y_1\dots w_nA_{t_n}y_n. $$ 
This is valid since  
\begin{align*} 
&\left\langle \delta_{w_1}*S*\dots *\delta_{w_n}*S*\Delta_k, g\right\rangle 
\\ 
&= \int S(y_1)\dots S(y_n)
\Delta_k(x)g\left(\left(\prod_{i=1}^n w_iy_i\right)x\right)\, dy\, dx 
\\ 
&= \int \prod_{i=1}^n \left[\psi(t_i)\ell(\rho^jr(y_i))\delta_{t_i}F(y_i)
\right]
\Delta_k(x)g\left(\left(\prod_{i=1}^n w_iy_i\right)x\right)
\, dy\, dx\,\bar{d}t, 
\end{align*} 
which will coincide with the integral in \eqref{2.9} after a change of 
variables.  
\par 
Let $DH(y,t)$ be the $n\times n$ matrix whose $i$th column vector is  
$\partial^L_{t_i}H(y,t)$: 
$$DH(y,t)=\left(\partial^L_{t_1}H(y,t), \dots , \partial^L_{t_n}H(y,t)\right).  $$
Then, \eqref{2.8} follows from the two estimates: 
\begin{multline} \label{2.10}
\left|\iiint \Delta_k(x)G_1(y,t)
g(H(y,t)x)\prod_{i=1}^n\left(\psi(t_i)\ell(t_i,y_i)F(y_i)\right)
\, dy\,\bar{d}t\,dx\right|
\\ 
\leq C(\log\rho)^{n}\rho^{\epsilon(k+c)/s'}\|\ell\|_{\Lambda_s^{\eta/s'}}^n
\|F\|_s^{n}, 
\end{multline}
\begin{multline} \label{2.11}
\left|\iiint \Delta_k(x)G_2(y,t)
 g(H(y,t)x)\prod_{i=1}^n\left(\psi(t_i)\ell(t_i,y_i)F(y_i)\right)
 \, dy\,\bar{d}t\,dx\right|
\\ 
\leq C(\log\rho)^{n}\rho^{\delta\epsilon(k+c)/s'}
\|\ell\|_{\Lambda_s^{\eta/s'}}^n\|F\|_1^{n}  
\end{multline}
with  
\begin{gather*}
G_1(y,t)=\zeta_1\left(\rho^{-n\epsilon k}\det(DH(y,t))\right), 
\\ 
G_2(y,t)=\zeta_2\left(\rho^{-n\epsilon k}\det(DH(y,t))\right),  
\end{gather*} 
where $\zeta_1$ is a function in 
$C_0^\infty(\Bbb R)$ such that $0\leq \zeta_1\leq 1$, 
$\supp(\zeta_1)\subset [-1,1]$, $\zeta_1(t)=1$ 
for $t\in [-1/2,1/2]$, $\zeta_2=1-\zeta_1$, and $\delta, \epsilon$ are small 
 positive numbers.  
\par 
Proof of \eqref{2.10}. 
Since  $\|\Delta_k\|_1\leq C$, 
$\int \psi(t_i)\,dt_i/t_i\leq C\log\rho$ and 
$$\int|F(y_i)\ell(t_i,y_i)|^s\psi(t_i)\,dy_i\,dt_i/t_i\leq C(\log\rho)
\|\ell\|_{d_s}^s\|F\|_s^s $$ 
(see the proof of Lemma 2), 
by H\"{o}lder's inequality, it suffices to show that 
\begin{equation}\label{2.12}
\int_{D_0^n}\chi_{[0,1]}\left(\rho^{-kn\epsilon}
\left|\det(DH(y,t))\right|\right)\, dy\leq C\rho^{\epsilon(k+c)} 
\end{equation}  
uniformly in $t\in [1,\rho^2]^n$ and $w_1, \dots w_n\in B(0,C\rho^2)$. 
By \eqref{1.1} and \eqref{1.4}  
$$ \partial^L_{t_i}H(y,t)=t_i^{-1}C[Q_i](A_{t_i}X(y_i)), 1\leq i\leq n-1,  
\quad  \partial^L_{t_n}H(y,t)=t_n^{-1}A_{t_n}X(y_n)$$ 
where $Q_i=\prod_{j=i+1}^n w_jA_{t_j}y_j$ for $1\leq i\leq n-1$. 
Fixing $y_2, \dots, y_n$ and changing variables with respect to $y_1$, we 
see that the integral in \eqref{2.12} is majorized by 
$$Ct_1^{-\gamma}\int_{\tilde{D}_0\times D_0^{n-1}}\chi_{[0,1]}
\left(c\rho^{-kn\epsilon}t_1^{-1}\left|\det(J(y,t))\right|
\right)\, dy_1d\bar{y},$$  
where $\tilde{D}_0=\{x\in \Bbb R^n: |x|\leq C\rho^{M}\}$, 
$M, C>0$, $d\bar{y}=dy_2\dots dy_n$ and 
$$J(y,t)=(y_1,\partial^L_{t_2}H(y,t), \dots , \partial^L_{t_n}H(y,t)).  $$ 
To see this, it may be convenient to write 
$$\partial^L_{t_1}H(y,t)=t_1^{-1}A_{\rho^2}C[A_{\rho^{-2}}Q_1](
A_{\rho^{-2}}A_{t_1}X(y_1)), $$ 
and to note that $|A_{\rho^{-2}}Q_1|\leq C$ (see \eqref{1.2}). 
Repeating this argument successively for  $y_2, \dots, y_n$, we can see that 
\eqref{2.12} follows from 
\begin{equation}\label{2.13}
(t_1\dots t_n)^{-\gamma}\int_{\tilde{D}_0^n}\chi_{[0,1]}
\left(c\rho^{-kn\epsilon}(t_1\dots t_n)^{-1}\left|\det Y\right|
\right)\, dy\leq C\rho^{\epsilon(k+c)},  
\end{equation}  
where   
$Y$ denotes the $n\times n$ matrix whose $i$th column vector is $y_i$. 
Write $y_i=(y_i^1, \dots, y_i^n)$.  
\par  
To prove \eqref{2.13}, expand $\det Y=\sum_{m=1}^n y_1^m
\Delta_{m1}$, where $\Delta_{m1}$ denotes the $(m,1)$ cofactor of $Y$. Then, 
using this and applying a rotation in $y_1$ variable, we see that 
\begin{align*}
&\int_{\tilde{D}_0^n}\chi_{[0,1]}\left(c\rho^{-kn\epsilon}
(t_1\dots t_n)^{-1}\left|\det Y 
\right|\right)\, dy 
\\
&= \int_{\tilde{D}_0^n}\chi_{[0,1]}
\left(c\rho^{-kn\epsilon}(t_1\dots t_n)^{-1}\left|y_1^1
\left(\sum_{m=1}^n\Delta_{m1}^2\right)^{1/2}\right|\right)\, dy 
\\ 
&\leq \int_{\tilde{D}_0^n}\chi_{[0,1]}
\left(c\rho^{-kn\epsilon}(t_1\dots t_n)^{-1}\left|y_1^1\Delta_{11}\right|\right)\, dy.  
\end{align*}  
Let $\tilde{D}_{01}=\{y_1\in \tilde{D}_0: |y_1^1|<t_1\rho^{k\epsilon}\}$, 
$\tilde{D}_{02}=\{y_1\in \tilde{D}_0: |y_1^1|\geq t_1\rho^{k\epsilon}\}$. 
Then we have 
\begin{align*} 
&\int_{\tilde{D}_0}\chi_{[0,1]}
\left(c\rho^{-kn\epsilon}(t_1\dots t_n)^{-1}\left|y_1^1\Delta_{11}\right|
\right)\, dy_1 
\\ 
&=\sum_{i=1}^2\int_{\tilde{D}_{0i}}\chi_{[0,1]}
\left(c\rho^{-kn\epsilon}(t_1\dots t_n)^{-1}\left|y_1^1\Delta_{11}\right|
\right)\, dy_1 
\\ 
&\leq Ct_1\rho^{k\epsilon}\rho^{M(n-1)}+ 
C\rho^{Mn}\chi_{[0,1]}\left(c\rho^{-k(n-1)\epsilon}(t_2\dots t_n)^{-1}
\left|\Delta_{11}\right|\right),   
\end{align*} 
and hence 
\begin{multline*}
\int_{\tilde{D}_0^n}\chi_{[0,1]}\left(c\rho^{-kn\epsilon}(t_1\dots t_n)^{-1}
\left|\det Y\right|\right)\, dy 
\\ 
\leq C\rho^{k\epsilon}\rho^b+ C\rho^{Mn}
\int_{\tilde{D}_0^{n-1}}\chi_{[0,1]}\left(c\rho^{-k(n-1)\epsilon}
(t_2\dots t_n)^{-1}\left|\Delta_{11}\right|\right)\, d\bar{y}  
\end{multline*}   
for some $b>0$.    
Repeating a procedure similar to this $n-1$ times,  we reach the estimate 
\begin{align*} 
&\int_{\tilde{D}_0^n}\chi_{[0,1]}\left(c\rho^{-kn\epsilon}(t_1\dots t_n)^{-1}
\left|\det Y\right|\right)\, dy 
\\ 
&\leq C\rho^{k\epsilon}\rho^b+ C\rho^b\int_{\tilde{D}_0}\chi_{[0,1]}
\left(c\rho^{-k\epsilon}t_n^{-1}\left|y_{n}^n\right|\right)\, dy_n   
\leq  C\rho^{k\epsilon}\rho^\tau 
\end{align*}  
for some $\tau>0$.  
This proves \eqref{2.13}.  
\par 
Proof of \eqref{2.11}.  
Let 
$$\tilde{\ell}(t_i,y_i)=\int_{s_i<t_i/2} 
\ell(t_i-s_i,y_i)\varphi_{\rho^{\epsilon k}}(s_i)\,ds_i, $$ 
where $\varphi_{u}(s_i)=u^{-1}
\varphi({u^{-1}}s_i)$, $u>0$,  with $\varphi\in C^\infty(\Bbb R)$ 
satisfying $\supp(\varphi)\subset (0,1/8)$, $\varphi\geq 0$, 
$\int \varphi(s)\,ds=1$.   Then 
\begin{gather*} 
\int\psi(t_i)|\ell(t_i,y_i)|\,dt_i/t_i\leq C(\log\rho)\|\ell\|_{d_1}, 
\\ 
\int\psi(t_i)|\tilde{\ell}(t_i,y_i)|\,dt_i/t_i\leq C(\log\rho)\|\ell\|_{d_1}, 
\\ 
\int\psi(t_i)|\ell(t_i,y_i)-\tilde{\ell}(t_i,y_i)|\,dt_i/t_i\leq C(\log\rho)
\omega(\ell,\rho^{\epsilon k}). 
\end{gather*}  
Therefore, writing 
\begin{align*}  
&\ell(t_1,y_1)\dots \ell(t_n,y_n)-\tilde{\ell}(t_1,y_1)\dots 
\tilde{\ell}(t_n,y_n) 
\\ 
&=(\ell(t_1,y_1)-\tilde{\ell}(t_1,y_1))\ell(t_2,y_2)\dots \ell(t_n,y_n)  
\\ 
&+ \tilde{\ell}(t_1,y_1)(\ell(t_2,y_2)-\tilde{\ell}(t_2,y_2))\ell(t_3,y_3)\dots \ell(t_n,y_n) 
\\ 
&+ \dots + \tilde{\ell}(t_1,y_1)\dots \tilde{\ell}(t_{n-1},y_{n-1})
(\ell(t_n,y_n)-\tilde{\ell}(t_n,y_n)),   
\end{align*} 
and applying the inequality $\omega(\ell,t)\leq \|\ell\|_{\Lambda_s^{\eta/s'}}
t^{\eta/s'}$, we see that 
to get  \eqref{2.11} it suffices to prove a variant of \eqref{2.11} where 
each $\ell(t_i,y_i)$ is replaced by $\tilde{\ell}(t_i,y_i)$ for $i= 1,2,\dots , n$. 
To show the estimate,  it suffices to prove  
\begin{multline}\label{2.14} 
\left|\iint \Delta_k(x)G_2(y,t)
 g(H(y,t)x)\prod_{i=1}^n\left(\psi(t_i)\tilde{\ell}(t_i,y_i)\right)
 \,\bar{d}t\,dx\right| 
 \\ 
\leq C\rho^{\delta\epsilon k}\rho^\tau\|\ell\|_{d_1}^n    
\end{multline} 
uniformly in $y\in D_0^n$ and $w_1, \dots, w_n\in B(0,C\rho^2)$ 
with some $\tau>0$, since the quantity on the left 
hand side of \eqref{2.14} is also 
bounded by $C(\log\rho)^n\|\ell\|_{d_1}^n$.  
\par 
To prove \eqref{2.14} we need the following three lemmas. 
\begin{lemma} 
Let $f$ be a continuous function  on $\Bbb R^n$ such that 
$$\supp(f)\subset B(0, C_1), \quad \int f(x)\,dx=0, \quad \|f\|_1\leq C_2. $$ 
Then there exist functions $f_1, f_2, \dots , f_n$ such that 
\begin{gather*}
f(x)=\sum_{i=1}^n\partial_{x_i}f_i(x), 
\\ 
\supp(f_i)\subset B(0, C_1'), \quad  \|f_i\|_1\leq C_2' \quad 
\text{for $i=1,2,\dots, n$,} 
\end{gather*}
with some constants $C_1'$ and $C_2'$.  
\end{lemma} 
This is from Lemma 7.1 in \cite{T}.  
\begin{lemma} 
Let $\Delta_k$ be as in \eqref{2.11}.  Then, 
there exist functions $F_j$, $j=1,2,\dots, n$, such that 
$\supp(F_j)\subset B(0,C\rho^k)$, $\|F_j\|_1\leq C\rho^{k\alpha}$ for some 
$\alpha>0$  and 
$$\Delta_k(x)=\sum_{j=1}^n\partial_{x_j} F_j(x).  $$ 
\end{lemma}  
This follows from Lemma 3. 
\begin{lemma} Suppose that $\det(DH(y,t)x)\neq 0$, where $DH(y,t)x$ is defined 
in the same way as $DH(y,t)$ with $H(y,t)x$ in place of $H(y,t)$.  Then   
$$\partial_{x_i}g(H(y,t)x)=\left\langle\nabla_t g(H(y,t)x), 
(DH(y,t)x)^{-1}(\partial_{x_i}^L(H(y,t)x))\right\rangle. $$ 
$($see Lemma $7.2$ of \cite{T}$)$.  
\end{lemma} 
\par 
By Lemma 4,  \eqref{2.14} follows from the estimate 
\begin{equation}\label{2.16} 
\left|\iint \partial_{x_m} F_m(x)g(H(y,t)x)a(t)
\prod_{i=1}^n\tilde{\ell}(t_i,y_i)\,\bar{d}t\,dx\right|
\leq C\rho^{k\delta\epsilon}\rho^\tau\|\ell\|_{d_1}^n   
\end{equation} 
for each $m$, $1\leq m\leq n$,  where
$a(t)=G_2(y,t)\prod_{i=1}^n\psi(t_i)$.  
Applying integration by parts and using the $L^1$ norm estimate for $F_m$ 
in Lemma 4, to prove \eqref{2.16} it suffices to show that 
\begin{equation}\label{2.17} 
\left|\int \partial_{x_m} g(H(y,t)x)a(t)\prod_{i=1}^n\tilde{\ell}(t_i,y_i)
\,\bar{d}t\right|
\leq C\rho^{-2nk\epsilon}\rho^\tau \|\ell\|_{d_1}^n
\end{equation} 
for all $x\in B(0,C\rho^k)$ with a sufficiently small $\epsilon>0$. 
By Lemma 5, the estimate \eqref{2.17} follows from 
\begin{multline}\label{2.18} 
\left|\int \left\langle\nabla_t g(H(y,t)x), 
(DH(y,t)x)^{-1}(\partial_{x_m}^L x)\right\rangle
a(t)\prod_{i=1}^n\tilde{\ell}(t_i,y_i)\,\bar{d}t\right| 
\\ 
\leq C\rho^{-2nk\epsilon}\rho^\tau\|\ell\|_{d_1}^n,   
\end{multline} 
since $\partial_{x_m}^L(H(y,t)x)=\partial_{x_m}^L x$ (see Section 2).  
 Note that 
$|\nabla_t a(t)|\leq C\rho^{-kn\epsilon}\rho^\tau$ and 
$|\det(DH(y,t)x)|\geq C\rho^{kn\epsilon}$ on the support of $a$, since 
$$|\det(DH(y,t)x)|=|\det \bar{C}[x]\det DH(y,t)|\geq C|\det DH(y,t)|,  $$  
where $\bar{C}[x]$ denotes the matrix expression for the linear transformation 
$C[x]$ (see \eqref{1.1}, \eqref{1.2}).  
Thus, taking into account Cramer's formula,  we have 
\begin{multline*}
\left|\partial_{t_u}\left[a(t)\prod_{i=1}^n\tilde{\ell}(t_i,y_i)
(DH(y,t)x)^{-1}(\partial_{x_m}^L x)\right]\right| 
\\ 
\leq C\rho^{-2nk\epsilon}\rho^\tau\prod_{i=1}^n|\tilde{\ell}(t_i,y_i)|+
C\rho^{-nk\epsilon}\rho^\tau|\partial_{t_u}\tilde{\ell}(t_u,y_u)|
\prod_{i\neq u}|\tilde{\ell}(t_i,y_i)|
\end{multline*}
for some $\tau>0$. Also, note that 
\begin{equation*}
\int_1^{\rho^2}|\tilde{\ell}(t_i,y_i)|\,dt_i/t_i
\leq C(\log\rho)\|\ell\|_{d_1}  
\end{equation*}
and   
$$\int_1^{\rho^2}|\partial_{t_i}\tilde{\ell}(t_i,y_i)|\,dt_i/t_i
\leq C(\log\rho)\rho^{-k\epsilon}\|\ell\|_{d_1},   $$ 
which follows from 
 $$|\partial_{t_i}\tilde{\ell}(t_i,y_i)|\leq C\rho^{-k\epsilon}
 \int |\ell(t_i-s_i,y_i)||\varphi'_{\rho^{\epsilon k}}(s_i)|\,ds_i. $$ 
These estimates along with  integration by parts imply \eqref{2.18}.   
This completes the proof of \eqref{2.11} and hence that of  Lemma 1. 
\end{proof} 

\section{Proof of Theorem 1}   
We use the following weighted Littlewood-Paley inequalities. 
\begin{lemma} Let $w\in \mathscr{A}_p$, $1<p<\infty$, and let the functions 
$\Delta_k$ be as in Section $3$.  Then 
\begin{gather}
\left\|\sum_k f_k*\Delta_k\right\|_{L^p(w)}\leq C_{p,w}
\left\|\left(\sum_k|f_k|^2\right)^{1/2}\right\|_{L^p(w)},    \label{3.+0}
\\ 
\left\|\left(\sum_k|f*\Delta_k|^2\right)^{1/2}\right\|_{L^p(w)}\leq 
C_{p,w}\|f\|_{L^p(w)},      \label{3.+1}
\end{gather} 
where the constant $C_{p,w}$ is independent of $\rho\geq 2$. 
\end{lemma}
\begin{proof} 
Let $K(x)=\sum_{k\in \mathscr I} \sigma_k\Delta_k(x)$, where $\mathscr I$ 
is an arbitrary finite subset of $\Bbb Z$ and $\{\sigma_k\}$ is an arbitrary 
sequence such that $\sigma_k= 1$ or $-1$.  Let $Sf(x)=f*K(x)$.  
Then 
\begin{enumerate}  
\renewcommand{\labelenumi}{(\arabic{enumi})}  
\item $S$ is bounded on $L^2$ with the operator norm bounded by 
a constant independent of $\rho$, $\mathscr I$ and $\{\sigma_k\}$;  
\item $|K(x)|\leq Cr(x)^{-\gamma}$;  
\item 
there are positive constants $C_1$ and $\epsilon$ such that $r(x)> C_1r(y)$ 
implies 
\begin{equation*} 
|K(y^{-1}x)-K(x)|\leq Cr(y)^\epsilon r(x)^{-\gamma-\epsilon}.    
\end{equation*}    
 \end{enumerate}    
The proof of \eqref{2.2} applies to show  
$$\|\Delta_k*\Delta_j\|_1\leq C\min\left(1,\rho^{-\epsilon|j-k|+c}\right). $$  
By this and the Cotlar-Knapp-Stein lemma we get (1). 
The estimates in (2) and (3) can be shown by a straightforward computation. 
We note that, to prove the estimate in (3), it suffices to show that 
if $r(x)>C_1r(y)$, then 
\begin{equation*} 
|\delta_k\phi(y^{-1}x)-\delta_k\phi(x)|\leq 
Cr(y)^\epsilon r(x)^{-\gamma-\epsilon}    
\end{equation*}    
for each $k$.   By application of dilation, this follows from the case $k=0$, 
which can be easily proved.   
\par 
Using (1), (2), (3) and applying  methods of \cite[Chapitre IV]{CM} and 
the proof of Theorem III in \cite{CF}, we have 
$\|Sf\|_{L^p(w)}\leq C\|f\|_{L^p(w)}$ for $w\in \mathscr{A}_p$, $1<p<\infty$,  
with a constant $C$ independent of $\mathscr I$, $\{\sigma_k\}$ and $\rho$. 
From this and the Khintchine inequality,  
\eqref{3.+1} follows.  A duality argument and \eqref{3.+1} imply 
\eqref{3.+0}.   
\end{proof} 
\par 
Let 
$$M_{F,\ell}f(x)=\sup_j|f*S_j(|F|,|\ell|)(x)|, $$ 
where $S_j(F,\ell)$  is as in Section 3 (see \eqref{3.a2}).  
Put $\mu^*f=M_{F,\ell}f$.  
Let $\theta\in (0,1)$.  
We prove the following result along with  Theorem 1. 
\begin{lemma} Let $s>1$,  $F\in L^s(D_0)$ and $\ell\in \Lambda_s^{\eta/s'}$ 
for some fixed $\eta>0$. 
Then, there exist positive constants 
$\epsilon$, $C$ independent of $\rho$ and $s$ such that 
$$\|\mu^*f\|_p\leq C(\log\rho)(1-\rho^{-\theta\epsilon/(2s')})^{-4/p}
\|\ell\|_{\Lambda_s^{\eta/s'}}\|F\|_s\|f\|_p  $$
for  $p> 1+\theta$.  
\end{lemma} 
In Lemmas 1 and 7, we can have the same value of $\epsilon$.   
\begin{proof}[Proof of Lemma $7$]  
Let  $U_\sigma =U_\sigma(F,\ell)$ (see \eqref{epsilon}) and write  
$U_\sigma f=\sum_{k_1,k_2} U_{k_1,k_2}f$, where 
$$U_{k_1,k_2}f
=\sum_j\sigma_jf*\Delta_{k_1+j}*\nu_j*\Delta_{k_2+j}, \quad 
\nu_j=\nu_j(F,\ell). $$ 
Fix integers $k_1, k_2$. 
By Lemma 1 and duality we have 
\begin{equation*}
\left\|f*\Delta_k*\nu_j\right\|_2\leq 
C(\log\rho)\min(1, \rho^{-\epsilon(|k-j|-c)/s'})
\|\ell\|_{\Lambda_s^{\eta/s'}}\|F\|_s\|f\|_2. 
\end{equation*} 
Using this along with Lemma 1,  for $\nu_j$ and $\tilde{\nu_j}$, 
and noting that  
 $\|\Delta_{k_2+j}*\Delta_{k_2+j'}\|_1\leq 
 C\min(1, \rho^{-\epsilon(|j-j'|-c)})$, 
 where we may assume that the number $\epsilon$ is equal to the value of 
 $\epsilon$ in Lemmas 1 and 7,  
 $\|\Delta_{k}\|_1\leq C$, $\|\nu_j\|_1\leq C(\log\rho)\|\ell\|_{d_1}\|F\|_1$, 
we have    
\begin{multline}\label{4.3} 
\left\|f*(\Delta_{k_1+j}*\nu_j)*(\Delta_{k_2+j}*\Delta_{k_2+j'})
*(\tilde{\nu}_{j'}*\Delta_{k_1+j'})\right\|_2
\\ 
\leq CA^2\min(1, \rho^{-2\epsilon(|k_1|-c)/s'}) 
\min(1, \rho^{-\epsilon(|j-j'|-c)})\|f\|_2,   
\end{multline} 
where $A=(\log\rho)\|\ell\|_{\Lambda_s^{\eta/s'}}\|F\|_s$,   
and also 
\begin{multline}\label{4.4} 
\left\|f*\Delta_{k_1+j}*(\nu_j*\Delta_{k_2+j})
*(\Delta_{k_2+j'}*\tilde{\nu}_{j'})*\Delta_{k_1+j'}\right\|_2
\\ 
\leq CA^2\min(1, \rho^{-2\epsilon(|k_2|-c)/s'})\|f\|_2.    
\end{multline} 
By \eqref{4.3} and \eqref{4.4}, taking the geometric mean we have   
\begin{multline*} 
\left\|f*\Delta_{k_1+j}*\nu_j*\Delta_{k_2+j}*\Delta_{k_2+j'}*\tilde{\nu}_{j'}
*\Delta_{k_1+j'}\right\|_2 
\\ 
\leq CA^2\prod_{i=1}^2\min(1, \rho^{-\epsilon(|k_i|-c)/s'}) 
\min(1, \rho^{-\epsilon(|j-j'|-c)/2})\|f\|_2.  
\end{multline*} 
We can obtain a similar estimate for 
$$\|f*\Delta_{k_2+j'}*\tilde{\nu}_{j'}*\Delta_{k_1+j'}*\Delta_{k_1+j}
*\nu_j*\Delta_{k_2+j}\|_2.  $$
Therefore, by the Cotlar-Knapp-Stein lemma we see that    
\begin{equation}\label{4.5} 
\left\|U_{k_1,k_2}f\right\|_2\leq 
CA\prod_{i=1}^2\min(1, \rho^{-\epsilon(|k_i|-c)/(2s')})\|f\|_2
\end{equation}  
uniformly in $\sigma$. By \eqref{4.5} we have 
\begin{equation}\label{4.6}   
\|U_\sigma f\|_2\leq \sum_{k_1,k_2}\|U_{k_1, k_2}f\|_2 
\leq CA(1-\rho^{-\epsilon/(2s')})^{-2}\|f\|_2
\leq CAB\|f\|_2,              
\end{equation}
 where $B=(1-\rho^{-\theta\epsilon/(2s')})^{-2}$. 
\par 
We define a sequence $\{p_j \}_1^{\infty}$ by $p_1 = 2$ and  
$1/p_{j+1} = 1/2  + (1 - \theta)/(2p_j)$ for $j \geq 1$.  
Then, $1/p_j=(1-a^j)/(1+\theta)$,  where $a=(1-\theta)/2$, so  
$\{p_j\}$ is decreasingly converges to $1 + \theta$. 
For $m \geq 1$ we show that   
\begin{equation}\label{4.7}   
\left \|U_{\sigma}f\right \|_{p_m} 
\leq C_m AB^{2/p_m}\left \|f\right \|_{p_m}
\end{equation} 
uniformly in $\sigma$, for all $F$ and $\ell$ satisfying the assumptions of 
Lemma 7. 
For $m=1$, this is a consequence of \eqref{4.6}.  Fix $m\geq 1$ and assume 
\eqref{4.7} for this $m$.  Then, using it for $U_\sigma(|F|,|\ell|)$ and 
applying the Khintchine inequality,  we see that  
\begin{equation}\label{4.8} 
\|g(f)\|_{p_m}\leq CAB^{2/p_m}\|f\|_{p_m},  
\end{equation}  
where 
$$g(f)=\left(\sum_j|f*\nu_j(|F|,|\ell|)|^2\right)^{1/2}$$  
(note that $\omega(|\ell|,t)\leq C\omega(\ell,t)$).   
Let $\nu^*(f)=\sup_j|f*|\nu_j||$ and 
$\Phi^*(f)=\sup_j|f*\Phi_j(|F|,|\ell|)|$, where $\nu_j=\nu_j(F,\ell)$ 
as above.  Note that 
\begin{gather*}
\nu^*(f)\leq \mu^*(|f|)+\Phi^*(|f|)\leq g(|f|)+2\Phi^*(|f|), 
\\ 
\Phi^*(|f|)\leq C(\log\rho)\|\ell\|_{d_1}\|F\|_1Mf.   
\end{gather*}  
These estimates and \eqref{4.8} 
along with the Hardy-Littlewood maximal theorem (see \cite{C, CW, GGKK}) imply 
\begin{equation}\label{4.9} 
\|\nu^*(f)\|_{p_m}\leq CAB^{2/p_m}\|f\|_{p_m}.       
\end{equation}  
Define $r_m$ by $1/r_m-1/2=1/(2p_m)$. Then by \eqref{4.9} and the estimate 
$\|\nu_j\|_1\leq CA$ we have the vector valued inequality (see \cite{DR} and 
also \cite{S, S2})   
\begin{equation}\label{4.10} 
\left\|\left(\sum |g_k*\nu_k|^2\right)^{1/2}\right\|_{r_m} 
\leq CAB^{1/p_m}
\left\|\left(\sum |g_k|^2\right)^{1/2}\right\|_{r_m}.     
\end{equation}
By the Littlewood-Paley theory (see Lemma 6) and \eqref{4.10} we have 
\begin{align}\label{4.11} 
\|U_{k_1,k_2}f\|_{r_m}&\leq C
\left\|\left(\sum_j |f*\Delta_{k_1+j}*\nu_j|^2\right)^{1/2}
\right\|_{r_m}
\\ 
&\leq CAB^{1/p_m}\left\|\left(\sum_j |f*\Delta_{k_1+j}|^2
\right)^{1/2}\right\|_{r_m}                \notag
\\ 
&\leq CAB^{1/p_m}\|f\|_{r_m}.              \notag 
\end{align}
Interpolating between \eqref{4.5} and \eqref{4.11}, since 
$1/p_{m+1}=(1-\theta)/r_m + \theta/2$,  we see that 
\begin{equation}\label{4.a1}  
\|U_{k_1,k_2}f\|_{p_{m+1}}\leq CAB^{(1-\theta)/p_m}
\prod_{i=1}^2\min(1, \rho^{-\theta\epsilon(|k_i|-c)/(2s')})
\|f\|_{p_{m+1}}.  
\end{equation}  
Thus 
\begin{align}
\|U_\sigma f\|_{p_{m+1}}&\leq \sum_{k_1,k_2}
\|U_{k_1,k_2}f\|_{p_{m+1}} \notag 
\\ 
&\leq CAB^{(1-\theta)/p_m}(1-\rho^{-\theta\epsilon/(2s')})^{-2}\|f\|_{p_{m+1}} 
\notag 
\\ 
&\leq CAB^{2/p_{m+1}}\|f\|_{p_{m+1}}.  \notag 
\end{align} 
This proves \eqref{4.7} for all $m$ by induction. For any $p\in (1+\theta,2]$ 
there exists a positive 
integer $j$ such that $p_{j+1}<p\leq p_j$. So, interpolating between the 
estimates \eqref{4.7} with $m=j$ and $m=j+1$, we have 
\begin{equation}\label{4.12}
\|U_\sigma f\|_{p} \leq CAB^{2/p}\|f\|_{p}. 
\end{equation} 
Let $g(f)$ be as in \eqref{4.8}.   
 The estimate \eqref{4.12} implies  
$\|g(f)\|_{p} \leq CAB^{2/p}\|f\|_{p}$ for $p\in (1+\theta, 2]$, from which 
Lemma 7 for $p\in (1+\theta, 2]$ follows, since $\mu^*(f)
\leq g(f)+\Phi^*(f)$. 
For $p>2$ Lemma 7 follows from interpolation between the estimate for $p=2$ of 
Lemma 7 and the estimate 
$$ \|\mu^*(f)\|_\infty\leq C(\log\rho)\|\ell\|_{d_1}\|F\|_1\|f\|_\infty. $$  
This completes the proof of Lemma 7. 
\end{proof} 
\par 
 Theorem 1 is an immediate consequence of the following result. 
\begin{lemma} Let the functions $h$, $\Omega$ be as in Theorem $1$ and put 
$\delta(p)= |1/p-1/p'|$.  
Suppose that $p\in (1+\theta, (1+\theta)/\theta)$. 
Let $A=(\log\rho)\|h\|_{\Lambda_s^{\eta/s'}}\|\Omega\|_s$ and let 
$B$ be as above\/$:$ 
$B=(1-\rho^{-\theta\epsilon/(2s')})^{-2}$.  
Then 
$$\|Tf\|_p\leq CAB^{1+\delta(p)}\|f\|_p,  $$ 
where the constant $C$ is independent of $s>1$, $\Omega$, $h$ and 
$\rho\geq 2$. 
\end{lemma} 
\begin{proof}     
Since $Tf=U_{\sigma^*}(K_0,h)(f)$, where $\sigma^*=\{\sigma_j\}$ with 
$\sigma_j=1$ for all $j$, by \eqref{4.12} we have 
$$\|Tf\|_p\leq CAB^{2/p}\|f\|_p \quad\text{for $p\in (1+\theta, 2]$.}$$ 
Now, a duality argument using a estimate similar to this one for 
$T^*f=U_{\sigma^*}(\tilde{K}_0,h)(f)$  will imply the conclusion for 
all $p\in (1+\theta, (1+\theta)/\theta)$.  
\end{proof} 
\par 
\begin{proof}[Proof of Theorem $1$] Take $\rho=2^{s'}$ in Lemma 8. Then 
$$\|Tf\|_p\leq Cs'(1-2^{-\theta\epsilon/2})^{-2(1+\delta(p))}
\|h\|_{\Lambda_s^{\eta/s'}}\|\Omega\|_s\|f\|_p $$ 
for $p\in (1+\theta, (1+\theta)/\theta)$ and $s>1$. 
Since $(1+\theta, (1+\theta)/\theta)
\to (1,\infty)$ as $\theta \to 0$, Theorem 1 follows from this estimate. 
\end{proof}

\section{Proof of Theorem 3} 
We need the following result to prove Theorem 3. 
\begin{lemma} Let $h$, $\Omega$ be as in Theorem $3$. 
Let $\theta\in (0,1)$ and 
 $A=(\log \rho)\|h\|_{\Lambda_s^{\eta/s'}}\|\Omega\|_{s}$.    We define  
\begin{equation*} 
R(f)(x) = \sup_{k \in \Bbb Z} \left|\sum_{j=k}^{\infty} 
f*S_jL(x) \right|,       
\end{equation*}  
where $S_jL$ is as in Section $3$. 
Let $I_\theta =(2(1+\theta)/(\theta^2-\theta+2), (1+\theta)/\theta)$.    
 Then, for   $p \in I_\theta$    we have  
$$\|R(f)\|_p \leq  CA
\left((1-\rho^{-\theta\delta/s'})^{-2(1 +\delta (p))}
+ (1-\rho^{-\theta\delta/s'})^{-4/p-1-\theta} \right)\|f\|_p $$ 
with some $\delta>0$,   
where   $C$ is independent of $s>1$, $h\in \Lambda_s^{\eta/s'}$, 
 $\Omega\in L^s(\Sigma)$ and $\rho$. 
\end{lemma}   
\begin{proof}      
 Let $\varphi_k=\sum_{m\geq k+2}\Delta_m=\delta_{\rho^{k+1}}\phi$.   
Using the decomposition  
$$\sum_{j=k}^{\infty}f*S_jL =T(f)*\varphi_k-
\left(\sum_{j=-\infty}^{k-1}f*S_jL\right)*\varphi_k + 
  \left(\sum_{j=k}^{\infty}f*S_jL\right)*(\delta-\varphi_k), $$  
we have  
\begin{equation}\label{6.2}   
R(f) 
\leq \sup_k\left|T(f)*\varphi_k\right|+
\sup_k\left|\left(\sum_{j=-\infty}^{k-1}f*S_jL\right)*\varphi_k\right| 
+ \sum_{j=0}^\infty N_j(f), 
\end{equation} 
where   
$N_j(f) = \sup_k \left|\left(f*S_{j+k}L\right)*(\delta - \varphi_k)\right|$. 
Lemma 8 and the Hardy-Littlewood maximal theorem imply  that      
\begin{equation}\label{6.3}  
\left\|\sup_k\left|T(f)*\varphi_k\right|\right\|_p 
\leq C\|M(Tf)\|_p 
\leq CA(1-\rho^{-\theta\epsilon/(2s')})^{-2(1+\delta(p))} \|f\|_p  
\end{equation} 
for $p\in (1+\theta, (1+\theta)/\theta)$.   
Also,  Lemma 7 and the Hardy-Littlewood maximal theorem imply that   
\begin{equation}\label{6.4}  
\|N_j(f)\|_u \leq CA(1-\rho^{-\theta\epsilon/(2s')})^{-4/u}\|f\|_u
\quad \text{for $u >1+\theta$. }
\end{equation} 
\par 
On the other hand, 
\begin{equation}\label{6.5}
N_j(f) \leq 
\left(\sum_k \left|f*S_{j+k}L*(\delta - \varphi_k)\right|^2\right)^{1/2}.  
\end{equation} 
Let 
$$V_\sigma f=\sum_k \sigma_k f*S_{j+k}L*(\delta - \varphi_k), $$ 
where $\sigma=\{\sigma_k\}$, $\sigma_k=1$ or $-1$.  
We prove  
\begin{equation}\label{6.a1} 
\|V_\sigma f\|_2 
\leq CA(1-\rho^{-\delta/s'})^{-3}\min(1, \rho^{-\delta(j-c)/s'}) \|f\|_2
\end{equation} 
for some $\delta, c>0$, uniformly in $\sigma$. 
Estimates \eqref{6.5} and \eqref{6.a1} with Khintchine's inequality imply 
 \begin{equation}\label{6.a13} 
\|N_j(f)\|_2 
\leq CA(1-\rho^{-\delta/s'})^{-3}\min(1, \rho^{-\delta(j-c)/s'})\|f\|_2.  
\end{equation} 
To prove \eqref{6.a1}, 
we apply an argument similar to the one used to prove \eqref{4.5}. 
We prove the estimates 
\begin{multline}\label{6.a6} 
\|f*S_{j+k}L*(\delta - \varphi_k)*(\delta - \varphi_{k'})*
S_{j+k'}\tilde{L}\|_2 
\\ 
\leq CA^2(1-\rho^{-\delta/s'})^{-2}\min(1, \rho^{-\delta(j-c)/s'})
\min(1, \rho^{-\delta(|k-k'|-c)/s'})\|f\|_2,  
\end{multline} 
\begin{multline}\label{6.a12} 
\|f*(\delta - \varphi_{k'})*
S_{j+k'}\tilde{L}*S_{j+k}L*(\delta - \varphi_k)\|_2 
\\ 
\leq CA^2(1-\rho^{-\delta/s'})^{-4}
\min(1, \rho^{-\delta(|k-k'|-c)/s'})\min(1, \rho^{-\delta(j-c)/s'})\|f\|_2. 
\end{multline} 
for some $\delta, c>0$, where $S_{j+k'}\tilde{L}=S_{j+k'}(\tilde{K}_0,h)$. 
By the Cotlar-Knapp-Stein lemma, the estimates \eqref{6.a6} and \eqref{6.a12} 
imply  \eqref{6.a1}. 
\par 
To prove \eqref{6.a6}, note that $\delta-\varphi_k=\sum_{m\leq k+1}\Delta_m$. 
Therefore, 
\begin{multline}\label{6.a2} 
\|f*S_{j+k}L*(\delta - \varphi_k)*(\delta - \varphi_{k'})*
S_{j+k'}\tilde{L}\|_2 
\\ 
\leq \sum_{m\leq k+1, m'\leq k'+1} 
\|f*S_{j+k}L*\Delta_m*\Delta_{m'}*S_{j+k'}\tilde{L}\|_2.   
\end{multline} 
By Lemma 1 we see that 
\begin{multline}\label{6.a3}
\|f*(S_{j+k}L*\Delta_m)*(\Delta_{m'}*S_{j+k'}\tilde{L})\|_2 
\\ 
\leq CA^2\min(1, \rho^{-\epsilon(|j+k-m|-c)/s'}) 
\min(1, \rho^{-\epsilon(|j+k'-m'|-c)/s'})\|f\|_2.  
\end{multline}  
Also, we have 
\begin{multline}\label{6.a4} 
\|f*S_{j+k}L*(\Delta_m*\Delta_{m'})*S_{j+k'}\tilde{L}\|_2  
\\ 
\leq CA^2\min(1, \rho^{-\epsilon(|m-m'|-c)})\|f\|_2.  
\end{multline}  
The estimates \eqref{6.a3} and \eqref{6.a4} imply 
\begin{multline}\label{6.a5}
\|f*S_{j+k}L*\Delta_m*\Delta_{m'}*S_{j+k'}\tilde{L}\|_2 
\\ 
\leq CA^2\min(1, \rho^{-\epsilon(|j+k-m|-c)/(2s')}) 
\min(1, \rho^{-\epsilon(|j+k'-m'|-c)/(2s')}) 
\\ 
\times \min(1, \rho^{-\epsilon(|m-m'|-c)/2})\|f\|_2.  
\end{multline}  
By \eqref{6.a2} and \eqref{6.a5} we have \eqref{6.a6}. 
\par 
Similarly, 
\begin{align}\label{6.a7} 
&\|f*(\delta - \varphi_{k'})*
S_{j+k'}\tilde{L}*S_{j+k}L*(\delta - \varphi_k)\|_2 
\\ 
&\leq \sum_{m\leq k+1, m'\leq k'+1} 
\|f*\Delta_{m'}*S_{j+k'}\tilde{L}*S_{j+k}L*\Delta_m\|_2   \notag 
\\ 
 &\leq \sum_{m\leq k+1, m'\leq k'+1} \sum_{\ell, \ell'}
\|f*\Delta_{m'}*S_{j+k'}\tilde{L}*\Delta_\ell*\Delta_{\ell'}*S_{j+k}L*
\Delta_m\|_2.                       \notag  
\end{align}   
(See \cite[p. 1555]{T} for the idea of interposing 
$\Delta_\ell*\Delta_{\ell'}$  in the convolution product.)   
By Lemma 1 we have 
\begin{multline}\label{6.a8} 
\|f*(\Delta_{m'}*S_{j+k'}\tilde{L})*(\Delta_\ell*\Delta_{\ell'})*(S_{j+k}L*
\Delta_m)\|_2
\\ 
\leq CA^2\min(1, \rho^{-\epsilon(|j+k'-m'|-c)/s'}) 
\min(1, \rho^{-\epsilon(|j+k-m|-c)/s'}) 
\\ 
\times \min(1, \rho^{-\epsilon(|\ell-\ell'|-c)})\|f\|_2.         
\end{multline} 
Also, 
\begin{multline}\label{6.a9} 
\|f*\Delta_{m'}*(S_{j+k'}\tilde{L}*\Delta_\ell)*(\Delta_{\ell'}*S_{j+k}L)*
\Delta_m\|_2
\\ 
\leq CA^2\min(1, \rho^{-\epsilon(|j+k'-\ell|-c)/s'}) 
\min(1, \rho^{-\epsilon(|j+k-\ell'|-c)/s'})\|f\|_2.         
\end{multline} 
By \eqref{6.a8} and \eqref{6.a9}, 
\begin{multline}\label{6.a10} 
\|f*\Delta_{m'}*S_{j+k'}\tilde{L}*\Delta_\ell*\Delta_{\ell'}*S_{j+k}L*
\Delta_m\|_2 
\\ 
\leq CA^2\min(1, \rho^{-\epsilon(|j+k'-m'|-c)/(2s')}) 
\min(1, \rho^{-\epsilon(|j+k-m|-c)/(2s')}) 
\\ 
\times \min(1, \rho^{-\epsilon(|\ell-\ell'|-c)/2})
\min(1, \rho^{-\epsilon(|j+k'-\ell|-c)/(2s')}) 
\min(1, \rho^{-\epsilon(|j+k-\ell'|-c)/(2s')})\|f\|_2.         
\end{multline} 
Summation with respect to $\ell$, $\ell'$ in \eqref{6.a10} implies 
\begin{multline}\label{6.a11} 
\|f*\Delta_{m'}*S_{j+k'}\tilde{L}*S_{j+k}L*
\Delta_m\|_2 
\\ 
\leq CA^2(1-\rho^{-\delta/s'})^{-2}
\min(1, \rho^{-\delta(|k-k'|-c)/s'}) 
\\ 
\times \min(1, \rho^{-\epsilon(|j+k'-m'|-c)/(2s')}) 
\min(1, \rho^{-\epsilon(|j+k-m|-c)/(2s')})\|f\|_2          
\end{multline} 
for some $\delta, c>0$.  
By \eqref{6.a7} and \eqref{6.a11} we obtain \eqref{6.a12}.  
\par 
For $p \in I_\theta$ 
we can  find  $u \in (1+\theta, 2(1+\theta)/\theta)$ 
such that $1/p = (1 - \theta)/u + \theta/2$,  so   
 an interpolation between \eqref{6.4} and \eqref{6.a13} implies that
\begin{equation}\label{6.8}  
\|N_j(f)\|_p \leq CA(1-\rho^{-\theta\delta/s'})^{-4(1-\theta)/u-3\theta} 
\min(1, \rho^{-\theta\delta(j-c)/s'}) 
\|f\|_p 
\end{equation}   
for some $\delta, c>0$. 
\par 
Also, we need the following result.  
\begin{lemma} There exist positive constants $C$, $C_1$ independent of 
$\rho$        such that 
$$\left|\left(\sum_{j=-\infty}^{k-1}S_jL\right)*\varphi_k(x) \right|
\leq C(\log \rho)\|h\|_{d_1}\|K_0\|_1
\rho^{-(k+1)\gamma}\chi_{[0,C_1]}(\rho^{-k-1}r(x)). $$  
\end{lemma} 
\begin{proof} 
Since $\int S_jL=0$, for $j\leq k-1$ we have 
$$S_jL*\varphi_k(x)=\rho^{-(k+1)\gamma}
\int\left(\phi(A_{\rho^{-k-1}}y^{-1}A_{\rho^{-k-1}}x)
-\phi(A_{\rho^{-k-1}}x)\right)S_jL(y)\,dy. $$ 
Also, since $\supp(S_jL)\subset \{r(x)\leq 2\rho^{j+2}\}$ and 
$\supp(\varphi_k)\subset \{r(x)\leq \rho^{k+1}\}$,  it follows that 
$\supp(S_jL*\varphi_k)\subset \{r(x)\leq C_1\rho^{k+1}\}$.  
Therefore 
\begin{align*}
\left|S_jL*\varphi_k(x)\right|&\leq 
C\rho^{-(k+1)\gamma}\chi_{[0,C_1]}(\rho^{-k-1}r(x))
\int |A_{\rho^{-k-1}}y| |S_jL(y)|\,dy 
\\ 
&\leq C\rho^{-(k+1)\gamma}\chi_{[0,C_1]}(\rho^{-k-1}r(x)) 
\rho^{(-k-1+j+2)/\beta_1}(\log \rho)\|h\|_{d_1}\|K_0\|_1.   
\end{align*} 
Thus summing over $j\leq k-1$, we get the 
conclusion.   
\end{proof} 
\par 
By Lemma 10   
\begin{equation}\label{6.9}  
\sup_k\left|f*\left(\sum_{j=-\infty}^{k-1}S_jL\right)*\varphi_k\right| 
\leq C(\log \rho)\|h\|_{d_1}\|K_0\|_1Mf.  
\end{equation} 
So, to estimate the maximal function on the left hand side of \eqref{6.9}, 
we can use the Hardy-Littlewood maximal theorem. 
\par 
By \eqref{6.2}, \eqref{6.3}, \eqref{6.8} and 
\eqref{6.9},  for $p \in I_\theta$  we have
\begin{equation*}
\|R(f)\|_p \leq CA\left((1-\rho^{-\theta\delta/s'})^{-2(1+\delta(p))} 
+(1-\rho^{-\theta\delta/s'})^{-4(1-\theta)/u-3\theta-1} \right) \|f\|_p   
\end{equation*} 
for some $\delta>0$.   
This implies the conclusion of Lemma 9, since 
 $4(1-\theta)/u +3\theta+1= 4/p+1+\theta$.   
\end{proof} 
\begin{proof}[Proof of Theorem $3$]      
Note that $T_*(f)\leq 2R(f)+CM_{K_0,h}(|f|)$.    
Therefore, Lemma 7 and Lemma 9 imply that  
 $$\|T_*(f)\|_p \leq C (\log\rho)
 \left(1-\rho^{-\theta\delta/s'}\right)^{-6}
 \|h\|_{\Lambda_s^{\eta/s'}}\|\Omega\|_s\|f\|_p $$
for $p \in I_\theta$ with some $\delta>0$.   
 Using this with $\rho=2^{s'}$ and noting that $I_\theta \to (1,\infty)$ as 
 $\theta\to 0$, 
 we can get the conclusion of Theorem 3.  
\end{proof}  

\section{Proof of Theorem 5}  
Let $M_{F,\ell}$ be as in Section 4. 
We prove Theorem 5 along with  the following result.  
\begin{proposition}  
Let $F\in L^q(D_0)$ and $\ell \in \Lambda_q^\eta$ for some $q>1$ and $\eta>0$. 
Let $1<p<\infty$.  Then, we have the following$:$  
\begin{enumerate}  
\renewcommand{\labelenumi}{(\arabic{enumi})}  
\item  if $q'\leq p<\infty$ and $w\in \mathscr{A}_{p/q'}$, the operator    
$M_{F,\ell}$ is bounded on $L^p(w);$  
\item   
the operator $M_{F,\ell}$ is bounded on $L^p(w^{1-p})$ if $1<p\leq q$ and 
$w\in \mathscr{A}_{p'/q'}$.  
\end{enumerate}  
\end{proposition} 
We use results of Sections 3, 4 and 5 with $\rho=2$.  We also write 
$\|f\|_{L^p(w)}=\|f\|_{p,w}$.  First, we prove results of Theorem 5 for $T$. 
\begin{proof}[Proof of Proposition $1 (1)$]
Since $\|S_j(|F|,|\ell|)\|_q\leq C2^{-j\gamma/q'}\|\ell\|_{d_q}\|F\|_q$,  
by the proof of Lemma 2,  and $\supp(S_j(|F|,|\ell|))\subset 
\{2^j\leq r(x)\leq 2^{j+3}\}$,  
H\"{o}lder's inequality implies that 
$$M_{F,\ell}(f)\leq C\|F\|_q\|\ell\|_{d_q}M_{q'}f, $$ 
where $M_{q'}f=\left(M(|f|^{q'})\right)^{1/q'}$.  
From this and the Hardy-Littlewood maximal theorem it follows that 
$$\|M_{F,\ell}(f)\|_{p,w}
\leq C\|F\|_q\|\ell\|_{d_q}\left\|M_{q'}f\right\|_{p,w}\leq C_{p,w}\|f\|_{p,w}
$$
if $q'<p$ and $w\in \mathscr{A}_{p/q'}$. 
\par 
Next, we handle the case $p=q'>1$. Let $w\in \mathscr{A}_1$. If $s>q'$, 
then $w\in \mathscr{A}_1\subset \mathscr{A}_{s/q'}$ and hence what we have 
already proved implies 
\begin{equation}\label{5.1} 
\|M_{F,\ell}(f)\|_{s,w} \leq C_{s,w}\|f\|_{s,w}.  
\end{equation} 
If $1<r<q'$, then by Lemma 7 
\begin{equation}\label{5.2}  
\|M_{F,\ell}(f)\|_{r} \leq C_{r}\|f\|_{r}.   
\end{equation} 
Interpolating with change of measure between \eqref{5.2} and \eqref{5.1} 
 with $w$ replaced by $w^{1+\tau}$ 
for sufficiently small $\tau>0$, we get 
$$\|M_{F,\ell}(f)\|_{q',w} \leq C_{q',w}\|f\|_{q',w}.  $$
This proves Proposition 1 (1). 
\end{proof}
\par  
\begin{remark} 
If $q'<p$ in Proposition $1 (1)$, then the assumption 
$\ell\in \Lambda^\eta$ is not needed.  Also, we can replace the assumption for 
$\ell$ of Proposition $1$ with the condition that there exists $\ell^*\in d_q$, $q>1$,  such that $|\ell|\leq \ell^*$ and $\ell^*\in \Lambda^\eta$ for some 
$\eta>0$, keeping  the conclusion unchanged,  since 
$M_{F,\ell}(f)\leq  M_{F,\ell^*}(|f|)$. In particular, if $\ell \in d_\infty$, 
we can take a constant function as $\ell^*$.  
\end{remark}
\begin{lemma} Let $B_jf(x)=f*\nu_j(x)$, where $\nu_j=\nu_j(F,\ell)$, 
$F\in L^1(D_0)$, $\ell\in d_1$ $($see \eqref{epsilon}$)$.    
Consider the inequality   
\begin{equation}\label{5.3} 
\left\|\left(\sum_{j=-\infty}^\infty|B_jf_j|^2\right)^{1/2}\right\|_{p,w}
\leq C_{p,w}\left\|\left(\sum_{j=-\infty}^\infty|f_j|^2\right)^{1/2}
\right\|_{p,w}. 
\end{equation} 
\begin{enumerate}  
\renewcommand{\labelenumi}{(\arabic{enumi})}  
\item Suppose that $F$ and $\ell$ are as in Proposition $1$.  Let $\delta\in 
(0,1)$. If \eqref{5.3} holds for some $p\in (1,\infty)$ and 
$w\in \mathscr{A}_p$, then $U_\sigma=U_\sigma(F,\ell)$ is bounded on 
$L^p(w^{1-\delta})$ uniformly in $\sigma$ $($see \eqref{epsilon}$)$.   
\item If $M_{F,\ell}$ is bounded on $L^p(w)$ for some $1<p\leq 2$ and 
$w\in \mathscr{A}_p$, then 
\eqref{5.3} holds with these $p$ and  $w$.
\end{enumerate}  
\end{lemma}  
\begin{proof} 
As in Section 4, we decompose $U_\sigma$ of (1) as 
$U_\sigma f=\sum_{k_1,k_2} U_{k_1,k_2}f$. 
By \eqref{5.3} and Lemma 6 we have 
\begin{align}\label{5.4}  
\|U_{k_1,k_2}f\|_{p,w}&\leq C\left\|\left(\sum_j|f*\Delta_{k_1+j}*\nu_j|^2
\right)^{1/2}\right\|_{p,w}    
\\ 
&\leq C\left\|\left(\sum_j|f*\Delta_{k_1+j}|^2\right)^{1/2}\right\|_{p,w}  
\notag 
\\ 
&\leq C\|f\|_{p,w}.    \notag  
\end{align} 
On the other hand, by the proof of Lemma 7 (see \eqref{4.a1}) and duality 
we have   
 \begin{equation}\label{5.5}  
 \|U_{k_1,k_2}f\|_{p} \leq C2^{-\epsilon(|k_1|+|k_2|)}\|f\|_p  
 \end{equation} 
 for some $\epsilon>0$. 
 Interpolating with change of measure between \eqref{5.5}  and  \eqref{5.4}, 
we see that 
 $$\|U_{k_1,k_2}f\|_{p,w^{1-\delta}}\leq C2^{-\delta\epsilon(|k_1|+|k_2|)}
\|f\|_{p,w^{1-\delta}}
 $$ 
 for all $\delta\in (0,1)$.  This implies that 
  $$\|U_\sigma f\|_{p,w^{1-\delta}}\leq \sum_{k_1,k_2}
\|U_{k_1,k_2}f\|_{p,w^{1-\delta}} \leq C\|f\|_{p,w^{1-\delta}}, 
  $$   
which proves part (1). 
  \par 
  Suppose that $M_{F,\ell}$ is bounded on $L^p(w)$ for $1<p\leq 2$. Then 
 \begin{equation}\label{5.6}  
\left\|\left(\sum_j|M_{F,\ell}f_j|^p\right)^{1/p}\right\|_{p,w} 
  \leq C\left\|\left(\sum_j|f_j|^p\right)^{1/p}\right\|_{p,w}.    
  \end{equation} 
  Also, we have 
 \begin{equation}\label{5.7}  
\left\| \sup_j|M_{F,\ell}f_j|\right\|_{p,w}\leq 
  C\left\| \sup_j|f_j|\right\|_{p,w}.     
  \end{equation} 
  Interpolating between \eqref{5.6} and \eqref{5.7},    
  $$\left\|\left(\sum_j|M_{F,\ell}f_j|^2\right)^{1/2}\right\|_{p,w} 
  \leq C\left\|\left(\sum_j|f_j|^2\right)^{1/2}\right\|_{p,w}.  $$  
Now, \eqref{5.3} follows from this estimate and a vector valued inequality for 
the Hardy-Littlewood maximal operator (see \cite[pp. 265--267]{GGKK}, 
\cite{RRT}).    This proves part (2). 
\end{proof} 
\par 
Let $q\geq 2$.  If $q'\leq p\leq 2$, $p>1$, by Proposition 1 (1) 
and Lemma 11, $U_\sigma$ is bounded on $L^p(w^{1-\delta})$  
for $w\in \mathscr{A}_{p/q'}$, where 
$U_\sigma=U_\sigma(F,\ell)$ and $F$, $\ell$ satisfy the assumptions of 
Proposition 1.  Replacing $w$ by $w^{1+\tau}$ for sufficiently small $\tau>0$ 
and taking $\delta$ suitably, we see that 
$U_\sigma$ is bounded on $L^p(w)$.    
This boundedness  also holds for $p\in (2, \infty)$ by the extrapolation 
theorem  of Rubio de Francia \cite{Ru}.   
If $1<p\leq q$, $w\in \mathscr{A}_{p'/q'}$, then this implies 
that $U_\sigma$ is bounded on $L^{p'}(w)$. Obviously, this is also valid for 
$U_\sigma^*=U_\sigma(\tilde{F},\ell)$. 
Therefore, by duality we can see that $U_\sigma$ is bounded on  
$L^{p}(w^{1-p})$. Let $\Omega$, $h$ be as in Theorem 5. 
By taking $F=K_0$, $\ell=h$, $\sigma_j=1$ for all $j$ 
in the definition of $U_\sigma$, now we can see that Theorem 5 holds 
for $T$  when $q\geq 2$.  
\par 
Also, from  a result of previous paragraph it follows that if $q\geq 2$, 
$1<p\leq q$,  $w\in \mathscr{A}_{p'/q'}$ and  $F$, $\ell$ are as in 
Proposition 1, then $M_{F,\ell}$ is bounded on $L^{p}(w^{1-p})$, 
since $M_{F,\ell}f\leq g(f)+CMf$ 
by the proof of Lemma 7 and the boundedness of $g$ follows from the uniform 
boundedness in $\sigma$ of 
$U_\sigma=U_\sigma(|F|,|\ell|)$, where  
\begin{equation*}
g(f)=\left(\sum_j|f*\nu_j(|F|,|\ell|)|^2\right)^{1/2}. 
\end{equation*}   
Here we recall that $\omega(|\ell|,t)\leq \omega(\ell,t)$.  
This proves Proposition 1 (2) for $q\geq 2$.  
\par 
It remains to prove Theorem 5 (for $T$) and Proposition 1 (2) when $1<q<2$.   
\begin{lemma}     
Let $1<q<2$, $2<p<\infty$. Let $F\in L^q(D_0)$, $\ell\in d_q$.   
If $M_{|\tilde{F}|^{2-q},|\ell|^{2-q}}$ is bounded on 
$L^{(p/2)'}(w^{-(p/2-1)^{-1}})$ and $w\in \mathscr{A}_p$, then      
 $$ 
\left\|\left(\sum_{j=-\infty}^\infty|B_jf_j|^2\right)^{1/2}\right\|_{p,w}
\leq C_{p,w}\left\|\left(\sum_{j=-\infty}^\infty|f_j|^2\right)^{1/2}
\right\|_{p,w},  
$$ 
where $B_j$ is defined as in Lemma $11$ by the functions $F$, $\ell$. 
\end{lemma}
\begin{proof} 
It suffices to prove the conclusion for $B_j'$  in place of 
$B_j$, where $B_j'f=f*S_j(F,\ell)$,  
on account of a vector valued inequality for the Hardy-Littlewood maximal 
operator. 
Take a non-negative function $u$ in $L^{(p/2)'}(w)$ with norm $1$ such that 
\begin{equation}\label{5.8} 
\left\|\left(\sum_{j}|B_j'f_j|^2\right)^{1/2}\right\|_{p,w}^2= 
\int \left(\sum_{j}|B_j'f_j|^2\right)u(x)w(x)\, dx. 
\end{equation} 
We see that 
\begin{equation}\label{5.1a}
|B_j'f(x)|^2\leq C\|\ell\|_{d_q}^q\|F\|_q^q 
(|f|^2*S_j(|F|^{2-q},|\ell|^{2-q}))(x).  
\end{equation}  
This can be proved as follows.  
First,  the Schwarz inequality implies that 
 \begin{multline*}
 |S_j(F,\ell)(x)|^2 
 \\ 
 \leq C\int_0^\infty \psi_j(t)|\ell(r(x))\delta_tF(x)|^q\, dt/t 
 \int_0^\infty \psi_j(t)|\ell(r(x))\delta_tF(x)|^{2-q}\, dt/t.  
\end{multline*}  
Therefore,  using 
$$\int\int_0^\infty \psi_j(t)|\ell(r(x))\delta_tF(x)|^q\, dt/t\, dx \leq 
C2^{j\gamma(1-q)}\|\ell\|_{d_q}^q\|F\|_q^q,    
$$ 
again by the Schwarz inequality, we have  
 \begin{multline*}
 |f*S_j(F,\ell)(x)|^2 
 \\ 
 \leq C2^{j\gamma(1-q)}\|\ell\|_{d_q}^q\|F\|_q^q
 \int\int_0^\infty \psi_j(t)|f(y)|^2|\ell(r(y^{-1}x))\delta_tF(y^{-1}x)|^{2-q}
 \, dt/t\, dy. 
\end{multline*} 
This implies \eqref{5.1a}. Therefore, 
the integral in \eqref{5.8} is majorized by 
\begin{equation}\label{5.9}  
C\|\ell\|_{d_q}^q\|F\|_q^q \int \left(\sum_{j}|f_j(y)|^2\right)
M_{|\tilde{F}|^{2-q},|\ell|^{2-q}}(uw)(y)\, dy. 
\end{equation} 
By H\"{o}lder's inequality, the integral in \eqref{5.9} is bounded by 
$$ \left\|\left(\sum_{j}|f_j|^2\right)^{1/2}
\right\|_{p,w}^2\left(\int \left|M_{|\tilde{F}|^{2-q},|\ell|^{2-q}}(uw)(y)
\right|^{(p/2)'}w^{-(2/p)(p/2)'}(y)\, dy\right)^{1/(p/2)'}. 
$$ 
Since $-(2/p)(p/2)'=-(p/2-1)^{-1}$, from the boundedness of 
$M_{|\tilde{F}|^{2-q},|\ell|^{2-q}}$ the last integral is majorized, 
up to a constant factor,  by   
$$\int \left|u(y)w(y)\right|^{(p/2)'}w^{-(2/p)(p/2)'}(y)\, dy
=\|u\|_{(p/2)',w}^{(p/2)'}=1. $$
Collecting results, we get the conclusion. 
\end{proof} 
\par 
Let $c_n=1-(1/2)^n$, $n=0, 1, 2,\dots$.  Suppose that 
$q^{-1}\in (c_n,c_{n+1}]$, $n\geq 1$. Put $r=q/(2-q)$. Then $(2r)^{-1}=
q^{-1}-1/2$, $r^{-1}\in (c_{n-1},c_n]$.  
For $n\geq 1$, we consider the following:  
\begin{assertion} 
Theorem $5$ for $T$ and Proposition $1$ $(2)$ hold when 
$q^{-1}\in (c_{n-1},c_n]$.   
\end{assertion} 
Assuming  Proposition $1$ $(2)$ when $q^{-1}\in (c_{n-1},c_n]$, 
we prove $A(n+1)$ ($n\geq 1$).  This will prove Theorem 5 for 
$T$ and Proposition 1 (2) when $1<q<2$, since we have already proved $A(1)$.  
\par 
 Suppose that $q^{-1}\in (c_n,c_{n+1}]$, $w\in 
\mathscr{A}_{p/q'}$, $q'\leq p<\infty$. 
Then, $(p/2)'\leq (q'/2)'=q/(2-q)=r$.  
Let $F$, $\ell$ satisfy the assumptions of Proposition 1.   
Since $r^{-1}\in (c_{n-1},c_n]$, $p/q'=
(p/2)/r'$, $-(p/2-1)^{-1}=1-(p/2)'$ and $|\ell|^{2-q}\in d_r$, 
$\omega(|\ell|^{2-q},t)
\leq C\omega(\ell,t)^{2-q}$, $|F|^{2-q}\in L^r(D_0)$,  by what we assume 
($A(n)$ for Proposition 1 (2)),  $M_{|\tilde{F}|^{2-q},|\ell|^{2-q}}$ is 
bounded on $L^{(p/2)'}(w^{-(p/2-1)^{-1}})$. By Lemmas 11 and 12, 
$U_\sigma=U_\sigma(F,\ell)$ is bounded on $L^p(w^{1-\delta})$. 
From this,  boundedness of $U_\sigma$ on $L^p(w)$ follows as before.  
This implies $A(n+1)$ for Theorem 5 (1) concerning $T$ 
as in the case when  $q\geq 2$. 
\par 
Suppose that $1<p\leq q$, $w\in \mathscr{A}_{p'/q'}$. Then, since $q'\leq p'$, 
by a result in the previous paragraph, $U_\sigma$ is bounded on 
$L^{p'}(w)$. We can see that the same is true for $U_\sigma^*$. 
By duality   $U_\sigma$ is bounded on $L^{p}(w^{1-p})$. 
This implies the boundedness on 
$L^{p}(w^{1-p})$ of $T$ and $M_{F,\ell}$ as in the case for $q\geq 2$.  
This finishes proving $A(n+1)$, and hence completes 
the proof of Theorem 5 for $T$ and Proposition 1.  
\par Next, we prove Theorem 5 for $T_*$. 
Let $\Omega$, $h$, $p, q, w$ be as in Theorem 5 (1). 
 By \eqref{6.2}, Lemma 10 and Theorem 5 for $T$, we have 
\begin{align}\label{5.10}
\|R(f)\|_{p,w}&\leq C\|M(Tf)\|_{p,w}+C\|Mf\|_{p,w}+
C\sum_{j=0}^\infty \|N_j(f)\|_{p,w} 
\\ 
&\leq C\|f\|_{p,w} + C\sum_{j=0}^\infty \|N_j(f)\|_{p,w}. \notag 
\end{align}
Since $N_j(f)\leq CMM_{K_0,h}(|f|)$, 
\begin{equation}\label{5.11}  
\|N_j(f)\|_{p,w}\leq C\|f\|_{p,w} 
\end{equation}  
by Proposition 1. By \eqref{6.4} and \eqref{6.a13}
\begin{equation}\label{5.12}  
\|N_j(f)\|_{p}\leq C2^{-\epsilon j}\|f\|_{p}   
\end{equation}
for some $\epsilon>0$.  
So, interpolating with change of measure between \eqref{5.12} and \eqref{5.11} 
with $w^{1+\tau}$ in place of $w$ for sufficiently  small $\tau>0$, we have 
\begin{equation}\label{5.13}  
\|N_j(f)\|_{p,w}\leq C2^{-\epsilon j}\|f\|_{p,w}   
\end{equation}
for some $\epsilon>0$.  
Since $T_*(f)\leq CR(f)+CM_{K_0,h}(|f|)$, 
by \eqref{5.10}, \eqref{5.13} and Proposition 1 
we have the $L^p(w)$ boundedness of $T_*$.  
This proves Theorem 5 (1). Theorem 5 (2) can be proved similarly. 

\section{Proof of Theorem 2} 
We give a proof of Theorem 2 by applying Theorem 1. 
Define 
$$E_m=\{\theta\in \Sigma: 2^{m-1}<|\Omega(\theta)|\leq 2^{m}\}$$ 
for $m=2,3,\dots$ and 
$$E_1=\{\theta\in \Sigma: |\Omega(\theta)|\leq 2\}. $$  
Let $\Omega_m=\Omega\chi_{E_m}-S(\Sigma)^{-1}\int_{E_m}\Omega\,dS$. 
Then $\int_\Sigma \Omega_m\,dS =0$, $\Omega=\sum_{m=1}^\infty\Omega_m$. 
\par 
Fix $p\in (1,\infty)$ and an appropriate function $f$ with $\|f\|_p\leq 1$.
Write $U(h,\Omega)=\|Tf\|_p$,   where $Tf=\mathop{\mathrm{p.v.}}f*L$, 
$L(x)=h(r(x))\Omega(x')r(x)^{-\gamma}$. 
Since $h\in\Lambda$, we can write $h=\sum_{k=1}^\infty a_kh_k$, where $\{a_k\}$  and $h_k$ are as in the definition of the space $\Lambda$.  Then, we 
decompose 
\begin{equation}\label{7.1}  
h\Omega=\sum_{m=1}^\infty\left(\sum_{k=m+1}^\infty a_kh_k\Omega_m + 
\sum_{k=1}^m a_kh_k\Omega_m\right).  
\end{equation}
By the subadditivity of $U$ and Theorem 1 we have    
\begin{align}\label{7.2}  
U\left(\sum_{k=m+1}^\infty a_kh_k,\Omega_m\right)&\leq 
\sum_{k=m+1}^\infty a_kU\left(h_k,\Omega_m\right) 
\\ 
&\leq C\sum_{k=m+1}^\infty ka_k\|h_k\|_{\Lambda_{1+1/k}^{1/(1+k)}}
\|\Omega_m\|_{1+1/k}                      \notag
\\ 
&\leq C\sum_{k=m+1}^\infty ka_k\|\Omega_m\|_{1+1/m}\leq C\|\Omega_m\|_{1+1/m}, 
\notag 
\end{align} 
since $\|\Omega_m\|_{1+1/k}\leq C\|\Omega_m\|_{1+1/m}$ if $k>m$.  
On the other hand, 
\begin{align}\label{7.3}  
U\left(\sum_{k=1}^m a_kh_k,\Omega_m\right)&\leq 
\sum_{k=1}^m a_kU\left(h_k,\Omega_m\right) 
\\ 
&\leq C\sum_{k=1}^m a_k m\|h_k\|_{\Lambda_{1+1/m}^{1/(1+m)}}
\|\Omega_m\|_{1+1/m}                      \notag
\\ 
&\leq C\sum_{k=1}^m a_k m\|\Omega_m\|_{1+1/m}
\leq Cm\|\Omega_m\|_{1+1/m}, 
\notag 
\end{align} 
since $\|h_k\|_{\Lambda_{1+1/m}^{1/(1+m)}}\leq 
C\|h_k\|_{\Lambda_{1+1/k}^{1/(1+k)}}\leq C$ if $k\leq m$. 
\par 
Note that 
\begin{equation*} 
\|\Omega_m\|_u\leq C2^m e_m^{1/u}, \quad 
1\leq u<\infty,     
\end{equation*} 
where $e_m=S(E_m)$.  
Using this 
and applying Young's inequality,  we see that 
\begin{align}\label{7.4} 
&\sum_{m\geq 1} m\|\Omega_m\|_{1+1/m}
\leq  C\sum_{m\geq 1} m2^m e_m^{m/(m+1)}  
 \\ 
 &\leq 2C\sum_{m\geq 1}(m/(m+1))m2^{m(1+1/m)}e_m  
+2C\sum_{m\geq 1}m2^{-m-1}/(m+1)  \notag 
\\ 
&\leq C\sum_{m\geq 1}m2^{m}e_m +C     
\leq  C\int_{\Sigma}|\Omega(\theta)|\log(2+|\Omega(\theta)|)\, dS(\theta) +C.
\notag 
\end{align}    
Theoerm 2 follows from  \eqref{7.1}, \eqref{7.2}, \eqref{7.3} and \eqref{7.4}. 
\begin{remark} 
Let $\mathscr M_a$, $a>0$, be the family of functions $h$ on 
$\Bbb R_+$ such that there exist a sequence $\{h_k\}_{k=1}^\infty$ of functions on $\Bbb R_+$ and a sequence $\{a_k\}_{k=1}^\infty$ of non-negative real 
numbers  satisfying 
$h=\sum_{k=1}^\infty a_kh_k$,  
$h_k\in d_{1+1/k}$,  $\|h_k\|_{d_{1+1/k}}\leq 1$ uniformly in $k\geq 1$  
and 
$\sum_{k=1}^\infty k^a a_k<\infty$.  Then, the space  $\mathscr M_1$ 
can be used to form kernels of singular integrals with a minimum size 
condition that allows us to get $L^p$ boundedness of singular integrals 
defined by the kernels from results of \cite{S, S2} (see \cite{S4}).   
\end{remark}


\begin{thebibliography}{99} 
\bibitem{C} A. P. Calder\'on,  
 {\it Inequalities for the maximal function relative to a metric},   
Studia Math.  {\bf 57} (1976),  297--306.    
\bibitem{CT} A. P. Calder\'on and A. Torchinsky,  
 {\it Parabolic maximal functions associated with a distribution},   
Advances in Math.  {\bf 16} (1975),  1--64.  
\bibitem{CZ} A.~P.~Calder\'on and A.~Zygmund, 
 {\it On singular integrals},  Amer.~J.~Math. {\bf 78}  
(1956),  289--309.     
\bibitem{CDF} Y. Chen, Y. Ding and D. Fan, {\it A parabolic singular integral 
operator with rough kernel}, J. Aust. Math. Soc. {\bf 84} (2008),  163--179. 
\bibitem{Ch} M. Christ,   
{\it Hilbert transforms along curves I.  Nilpotent groups},  
 Ann. of Math. {\bf 122} (1985),   575--596.   
\bibitem{CF}  R. R. Coifman and C. Fefferman, 
{\it Weighted norm inequalities for maximal functions and singular 
integrals}, Studia Math. {\bf 51} (1974), 241--250.  
\bibitem{CM}  
 R. R. Coifman and Y. Meyer, 
 {\it Au del\`a des op\'erateurs pseudo-diff\'erentiels},  
  Ast\'erisque no. 57,   Soc. Math. France,  1978.   
\bibitem{CW}  R. R. Coifman and G. Weiss, 
 {\it Analyse Harmonique Non-Commutative sur Certains Espaces Homogenes},  
 Lecture Notes in Math. 242, Springer-Verlag, Berlin and New York, 1971.  
\bibitem{D} J.~Duoandikoetxea,  
{\it Weighted norm inequalities for homogeneous singular integrals},  
 Trans. Amer. Math. Soc. {\bf 336} (1993),  869--880.   
\bibitem{DR} 
 J.~Duoandikoetxea and J.~L.~ Rubio de Francia, 
 {\it Maximal and singular integral operators via Fourier transform 
 estimates},   Invent. Math. {\bf 84} (1986), 541--561. 
\bibitem{FR} E. B. Fabes and  N. Rivi\`{e}re,  
{\it Singular integrals with mixed homogeneity},  
 Studia Math. {\bf 27}  (1966),  19--38.   
\bibitem{FS} D. Fan and S. Sato,  
{\it Weighted weak type $(1,1)$ estimates for  singular
integrals and Littlewood-Paley functions},  
 Studia Math. {\bf 163}  (2004),  119--136.   
 \bibitem{GGKK} I. Genebashvili, A. Gogatishvili, V. Kokilashvili, 
 and M. Krbec, {\it Weight theory for integral transforms on spaces of 
 homogeneous type}, Pitman Monographs and Surveys in Pure and Appl. Math. 92, 
 Addison Wesley Longman, 1998. 
\bibitem{NRW} A. Nagel, N. Rivi\`{e}re and S. Wainger, {\it On Hilbert 
transforms along curves, II}, Amer. J. Math. {\bf 98} (1976), 395--403. 
\bibitem{NS} A. Nagel and E. M. Stein, {\it Lectures on pseudo-differential 
operators}, Mathematical Notes 24, Princeton University Press,  Princeton, NJ, 
1979.  
\bibitem{RS}  F. Ricci and E. M. Stein,  
{\it Harmonic analysis on nilpotent groups and singular integrals, I. 
Oscillatory integrals}, J. Func. Anal.  {\bf 73} (1987),  179--194.  
\bibitem{RS2}  F. Ricci and E. M. Stein,  
{\it Harmonic analysis on nilpotent groups and singular integrals, II. 
Singular kernels supported on submanifolds}, 
J. Func. Anal.  {\bf 78} (1988),  56--84.  
\bibitem{R} N. Rivi\`{e}re, 
{\it Singular integrals and multiplier operators},  
 Ark. Mat.   {\bf 9} (1971),    243--278.  
\bibitem{Ru} J.L. Rubio de Francia, 
{Factorization theory and $A_p$ weights},  
Amer. J. Math.  {\bf 106} (1984),  533--547.   
\bibitem{RRT} 
 J. L. Rubio de Francia, F. J. Ruiz and J. L. Torrea,  
{\it Calder\'on-Zygmund theory for operator-valued kernels},  
 Adv. in Math. {\bf 62} (1986), 7--48.  
 \bibitem{S} S.~Sato, 
{\it Estimates for singular integrals and extrapolation},   
Studia Math. {\bf 192} (2009), 219--233.  
  \bibitem{S2} S. Sato, 
{\it Estimates for  singular integrals along surfaces of revolution},  
J. Aust. Math. Soc. {\bf 86} (2009), 
413--430. 
\bibitem{S3} S. Sato, 
{\it Weak type $(1,1)$ estimates for parabolic  singular integrals}, 
 Proc. Edinb. Math. Soc. (to appear).  
\bibitem{S4} S. Sato, 
{\it A note on $L^p$ estimates for singular integrals}, Sci. Math. 
Jpn. {\bf 71} (2010), 343--348.  
\bibitem{Se}  A. Seeger,   
{\it Singular integral operators with rough convolution kernels}, 
 J.  Amer. Math. Soc.   {\bf  9} (1996),  95--105.  
 \bibitem{St}  E. M. Stein, 
 {\it Harmonic Analysis$:$  Real-Variable Methods, Orthogonality and 
Oscillatory Integrals},    Princeton University Press,   Princeton, NJ, 1993.  
\bibitem{SW}  E.~M.~Stein and S.~Wainger,  
 {\it Problems in harmonic analysis related to curvature},   
 Bull.~Amer.~Math.~Soc.  {\bf 84}  (1978),    1239--1295. 
\bibitem{T} T. Tao, 
{\it The weak-type $(1,1)$ of $L\log L$ homogeneous convolution operator},  
Indiana Univ. Math. J.  {\bf 48} (1999),  1547--1584.    
 \bibitem{W} D.~Watson, 
{\it Weighted estimates for singular integrals via Fourier transform 
estimates}, 
 Duke Math. J. {\bf 60} (1990),   389--399. 


\end{thebibliography}
\end{document}